\documentclass[twocolumn]{autart} 

\usepackage{graphicx}          

\usepackage[cmex10]{amsmath}
\usepackage{amsfonts}
\usepackage{amssymb}
\usepackage{verbatim} 
\usepackage{graphics,psfrag,epsfig}
\usepackage{epstopdf} 
\usepackage[noadjust]{cite}
\usepackage{cases}

\usepackage{color}
\definecolor{red}{rgb}{1,0.2,0.2}
\definecolor{green}{rgb}{0.2,1,0.5}
\definecolor{blue}{rgb}{0,0,1}
\definecolor{lightblue}{rgb}{0.3,0.5,1}

\newcommand{\diam}{\mathrm{diam}}
\newcommand{\dist}{\mathrm{dist}}

\newcommand{\0}{\mathbf{0}}
\newcommand{\1}{\mathbf{1}}

\newcommand{\bb}{\mathbf{b}}
\newcommand{\co}{\mathrm{conv}}

\newcommand{\g}{\mathbf{g}}

\newcommand{\p}{\mathbf{p}}
\newcommand{\piB}{\boldsymbol{\pi}}

\newcommand{\Rm}{\mathbb{R}^m}
\newcommand{\sB}{\mathbf{s}}
\newcommand{\T}{^{\mathsf{T}}}
\newcommand{\uu}{\mathbf{u}}
\newcommand{\vv}{\mathbf{v}}

\newcommand{\x}{{\mathbf{x}}}
\newcommand{\y}{\mathbf{y}}
\newcommand{\z}{\mathbf{z}}

\newcommand{\nnb}{\nonumber}

\newtheorem{theorem}{Theorem}
\newtheorem{lemma}{Lemma}
\newtheorem{proposition}{Proposition}
\newtheorem{corollary}{Corollary}
\newtheorem{assumption}{Assumption}
\newtheorem{remark}{Remark}

\begin{document}
	
	\begin{frontmatter}
		
		\title{Distributed Optimization over Directed Graphs with Row Stochasticity and Constraint Regularity\thanksref{footnoteinfo}} 
		
		\thanks[footnoteinfo]{This work
			was  supported  in  part  by  the Air Force Office of Scientific Research through MURI AFOSR Grant \#FA9550-09-1-0538. The  material in this paper was partially presented at the 2016 American Control Conference, July 6--8, Boston, MA, USA. 
			Corresponding author: E.~H.~Abed. Tel. +1 301 405 3631. 
			Fax +1 425 648 3960.}
		
		\author{Van Sy Mai}\ead{vsmai@terpmail.umd.edu},    
		\author{Eyad H. Abed}\ead{abed@umd.edu}               
		
		\address{University of Maryland, College Park, MD, 20742, USA}  

		\date{}
		
		\begin{abstract}                          
			This paper deals with an optimization problem over a network of agents, where the cost function is the sum of the individual (possibly nonsmooth) objectives of the agents and the constraint set is the intersection of local constraints. 
			Most existing methods employing subgradient and consensus steps for solving this problem require the weight matrix associated with the network to be column stochastic or even doubly stochastic, conditions that can be hard to arrange in directed networks. 
			Moreover, known convergence analyses for distributed subgradient methods vary depending on whether the problem is unconstrained or constrained, and whether the local constraint sets are identical or nonidentical and compact. 
			The main goals of this paper are: (i) removing the common column stochasticity requirement; (ii) relaxing the compactness assumption, and (iii) providing a unified  convergence analysis.
			Specifically, assuming the communication graph to be fixed and strongly connected and the weight matrix to (only) be row stochastic,  a distributed projected subgradient algorithm and a variation of this algorithm are presented to solve the problem for cost functions that are convex and Lipschitz continuous. 
			The key component of the algorithms is to adjust the subgradient of each agent by an estimate of its corresponding entry in the normalized left Perron eigenvector of the weight matrix. 
			These estimates are obtained locally from an augmented consensus iteration using the same row stochastic weight matrix and requiring very limited global information about the network. 
			Moreover, based on a regularity assumption on the local constraint sets, a unified analysis is given that can be applied to both unconstrained and constrained problems and without assuming compactness of the constraint sets or an interior point in their intersection. 
			Further, we also establish an upper bound on the absolute objective error evaluated at each agent's available local estimate under a nonincreasing step size sequence. This bound allows us to analyze the convergence rate of both algorithms. 
		\end{abstract}
		
	\end{frontmatter}
	
	\section{Introduction}\label{secIntroduction}
	We consider a network of agents without a central coordination unit that is tasked with solving a global optimization problem in which the objective function is the sum of local costs of the agents, that is, 
	$F(\x) = \sum_{i=1}^n f_i(\x)$ 
	where $f_i:\mathbb{R}^m \to \mathbb{R}$ represents the private objective of agent $i$ and $n$ is the number of agents in the network. In addition, each agent is  associated with a private constraint set. 
	This problem arises in many applications, such as distributed multi-agent coordination, estimation in sensor networks, resource allocation, and large-scale machine learning, to name a few. 
	Many distributed optimization methods have been developed to address this problem; see e.g., \cite{Tsitsiklis86, Bertsekas89, Rabbat04, Johansson08CDC, Nedic09AC, WE10, Boyd11, Lobel11AC, Duchi12AC, Nedic15AC, Lin16distributed, Olshevsky16, Makhdoumi17AC, Bajovic17} and references therein.

	Although much research has been carried out in this problem area, much of the existing literature invokes the assumption that communication among agents is bidirectional; see, e.g., \cite{Boyd11,Jakovetic15AC,Makhdoumi17AC} and references therein for classes of distributed alternating direction and augmented Lagrangian methods,  and \cite{Wang11,Shi15,Lee16AC, Olshevsky16} for (sub)gradient and consensus based algorithms.  
	However, this assumption is not always practical in real-world distributed systems where communications can be unidirectional, especially when a broadcast-based communication scheme is used; e.g., agents in a wireless network may have different communication ranges due to environmental effects or individual  broadcast power limits. 
	Motivated by this limitation, many recent works have considered directed communication by using weaker assumptions on how the agents diffuse and aggregate local information--a process which can be represented through a weight matrix that is in accordance with the network structure. 
	In particular, a common assumption in this line of research is the availability of a \emph{doubly stochastic} weight matrix (i.e., each row and each column sum up to 1)
	or a \emph{column stochastic} one (i.e., each column sums up to 1); see, e.g.,  \cite{Johansson08CDC,Nedic10AC,Ram10,Lobel11AC,Lobel11distributed,Tsianos12ACC,Chen12,Lee13,Nedich16arXiv} for the former case and 
	\cite{Tsianos12ACC, Tsianos12Allerton,Nedic15AC,Akbari15AC,Zeng15arXiv, Xi17AC} for the latter. 
	Here, the latter case is clearly  weaker and, more importantly, allows one to employ the Push-Sum protocol  \cite{Kempe03gossip} or the technique in \cite{Makhdoumi15} in order to asymptotically ``balance the graph," thereby achieving exact convergence as a doubly stochastic matrix does in many distributed algorithms. 
	However, the implementation of a column stochastic matrix  requires that each agent knows its out-degree exactly and/or controls its outgoing weights  (i.e., the weights that other agents put on its information) so that they sum up to one.  
	Such  a requirement, however, is 
	impractical in many situations, especially when agents use  broadcast-based communications and thus they neither know their out-neighbors nor are able to adjust their outgoing weights; a wireless \textit{ad hoc} network is such an example%
	
	
	In comparison with a column stochastic weight matrix, one that is \textit{row stochastic} (i.e., each row sums up to 1), is much easier to achieve in a distributed setting. Here, each agent can individually (and to some extent arbitrarily)  decide the weights on information it receives from its neighbors. Thus, if the weight matrix is required to be only row stochastic, there is no need for nodes to send acknowledgment signals. 
	This makes row stochastic matrices particularly suitable for reaching consensus in  broadcast-based communication environments. 
	However, when a row stochastic matrix is used for distributed optimization, most (sub)gradient based algorithms fail to achieve an optimal solution due to the nonuniform stationary distribution of the weight matrix (aka the normalized left Perron eigenvector). In \cite{Tsianos12ACC}, the authors suggest a re-weighting technique that makes it possible to use a row stochastic matrix in distributed optimization. The same technique is also employed in \cite{Lin16distributed}. 
	However, implementation of the algorithms in \cite{Tsianos12ACC,Lin16distributed} assumes \emph{global knowledge} of the graph, namely, the left Perron eigenvector and the number of agents in the network. 
	Moreover, when there is any change in the network structure, 
	the reweighting technique requires the whole network to be reprogrammed with a new Perron eigenvector, which may not be immediately available.
	
	In this connection, one may argue that one can employ available distributed algorithms for estimating this vector, e.g., \cite{QLL11, PGMS13,Charalambous16}, then apply the re-weighting technique. However, we note that the finite-time algorithm in \cite{Charalambous16}, relying on each agent keeping track of the rank of a  Hankel matrix that is growing in size, does not work for all  initial conditions and becomes inexact and less reliable as the network size increases (see \cite{Charalambous16} for a detailed discussion and numerical illustrations); not to mention the memory requirement and computational burden imposed on each agent. The algorithm in \cite{QLL11, PGMS13} is exact but converges asymptotically. Thus, it must be terminated after a finite number of steps prior to implementing any subsequent optimization algorithm, resulting in inexact convergence. Moreover, the needed number of steps depends on the desired estimation accuracy, as well as network size and structure, and thus needs to be determined (or redetermined if any change in the network structure occurs) by a central coordinator and made available to the agents before running the algorithm. Further, upon termination, only approximations of the eigenvector are available, which could still cause significant errors when used in the re-weighting technique. 
	Thus, a fully distributed algorithm employing only row stochastic weight matrices has not been available in the field of distributed optimization thus far. 
	
	

	In this work, we achieve such algorithms under mild requirements on available global network information.  
	More precisely, we present a distributed algorithm and a variation on the algorithm that use a row stochastic weight matrix and assume that the agents have different IDs and know only an upper bound on the network size. Our idea is as follows. 
	We let all the agents perform an augmented consensus protocol similar to that in \cite{PGMS13} in order to estimate the stationary distribution of the weight matrix while updating their states using an iteration akin to that in the Distributed Projected Subgradient (DPS) method (see, e.g., \cite{Nedic10AC, Ram10, Lee13}), except that subgradient values are now scaled appropriately and locally by the agents. 
	Here, the estimation step is implemented concurrently with the optimization step, and thus no communication overhead is added. 
	Moreover, although the algorithm is based on the DPS method, we believe that its principle (i.e., the use of a particular augmented consensus) can be generalized to a class of distributed algorithms that use consensus and subgradient.  
	
	Another important contribution is our unified convergence analysis (including the  convergence rate) of a distributed  subgradient method that applies to both  unconstrained and constrained problems, where the local convex cost functions may be non-differentiable and the private constraint sets satisfy certain regularity assumptions, so that they need not be identical or bounded or have an interior point in their intersection. 
	Existing analyses of subgradient-based methods differ according to whether the problem is unconstrained (see, e.g., \cite{Nedic09AC,Lobel11AC,Chen12,Nedic15AC,Olshevsky16}) or constrained with identical  constraint sets (see, e.g., \cite{Johansson08CDC,Ram10,Lobel11distributed,Duchi12AC,Tsianos12ACC,Mai16ACC2}); compactness also is needed in \cite[Thm.~5.2]{Ram10} and \cite{Mai16ACC2}. 
	Nonidentical constraints are considered in \cite{Nedic10AC,Lobel11distributed,Lin16distributed}, where the local constraint sets are assumed to be compact and their intersection has a nonempty interior. 
	Recent efforts have been made to deal with the case of nonidentical unbounded constraint sets. Specifically, \cite{Lee16AC} assumes a set-regularity condition  (which as shown later is stronger than ours)  and  differentiability of $f_i$ together with boundedness and Lipschitz continuity of its gradient. The work \cite{LinRen17arXiv}, relying on the nonempty interior assumption and differentiability of $f_i$, 
	proposes to use state-dependent step sizes so as to confine the agents' states to a bounded region.%
	
	Our third contribution is the derivation of the convergence rate of our  algorithms  for the case of nonidentical constraint sets  with possibly an empty-interior intersection. This, to the best of our knowledge, has not been established previously for subgradient methods in the same setting. The main challenge stems from the fact that the agents' local estimates need not be feasible at all times. 
	The paper \cite{Lee16AC} establishes only an error bound on the (expected) distance from the agents' iterative sequence to the optimal point under a constant step size, hence  inexact convergence. 
	Here, we provide a bound on the absolute objective error which demonstrates  how the rate of convergence depends on the step size sequence, exhibiting similarity to that of the centralized subgradient approach. 
	Specifically, a sublinear rate of $O({\ln (t)}/{\sqrt{t}} )$ is achieved when using a step size diminishing at rate $O(1/\sqrt{t})$, which is  standard in non-smooth optimization (see, e.g., \cite{Nesterov04,Chen12,Duchi12AC,Nedic15AC} and references therein). 
	Though better convergence rates are attainable when certain assumptions on smoothness and strong convexity of $f_i$ are imposed  \cite{WE10,Shi15,Nedich16arXiv,Xi17AC}, fast convergence is not the main goal in this paper; instead, we focus on a distributed  non-smooth optimization problem with row stochasticity and constraint regularity.


	Preliminary work along the lines of this paper appeared in \cite{Mai16ACC2}, where only one algorithm was presented and several proofs were omitted. 
	In addition, it is assumed in \cite{Mai16ACC2} that the  sets $X_i$ are identical and compact, while in this paper we consider nonidentical constraint sets and relax the compactness requirement, allowing for broader applicability. The present paper further introduces a variation on the algorithm presented in \cite{Mai16ACC2}, and presents a new convergence analysis that holds for both algorithms under these relaxations. Here the proof technique relies on  the regularity assumption on the local constraint set, differing significantly from that in \cite{Mai16ACC2}. 
	Finally, the convergence rate, which was not shown in \cite{Mai16ACC2}, is studied here for both algorithms. 
	
	
	The rest of the paper proceeds as follows. The problem formulation and proposed algorithms are given in Section \ref{secProblem}. 
	Convergence and convergence rate of the algorithms are studied in Sections \ref{secMainResults} and \ref{secRoC}, respectively. 
	Section \ref{secNumericalExample} includes a numerical example to illustrate the paper's findings. Conclusions are given in Section~\ref{secConclusion}. 
	
	\textit{Notation and basic terminology:} 
	Vectors are denoted in bold fonts, e.g., $\x = [x_1, x_2, ..., x_m]\T\in \mathbb{R}^m$, $\x_i = [x_{i1}, x_{i2}, ..., x_{im}]\T \in \mathbb{R}^m$, $\1 = [1,1,...,1]\T$ and $\e_i=[0,...,0,1_i,0,...,0]\T$. 
	For a vector $\x$, $\|\x\|$ denotes its 2-norm. 
	For a matrix $A$, $A\T$ denotes its transpose, $[A]_{ij}$ (also $a_{ij}$) the $(ij)$-th element, and $\|A\|$ the (induced) 2-norm of $A$. 
	A nonnegative square matrix $A$ is row stochastic if $A\1=\1$, column stochastic if $A\T\1=\1$,  doubly stochastic if it is both row and column stochastic. 
	
	A directed graph $\mathcal{G} {=} (\mathcal{V}, \mathcal{E})$ consists of a finite set of nodes $\mathcal{V} = \{1,2,...,n\}$ and a set $\mathcal{E} \subseteq \mathcal{V}\times \mathcal{V}$ of edges, where an ordered pair $(i,j) \in \mathcal{E}$ indicates that node $i$ receives information from node $j$. 
	A directed path is a sequence of edges of the form $(i_1, i_2), (i_2, i_3),..., (i_{k-1}, i_k)$. 
	Node $i$ is said to be reachable from node $j$ if there exists a path from $j$ to $i$. 
	Each node is reachable from itself. A graph 
	$\mathcal{G}$ is strongly connected if any node is reachable from any other node. See \cite{Godsil01} for more background on graphs. 
	%
	
	Let $f:\!\mathbb{R}^m \!\to \mathbb{R}$ be a convex function. 
	Then $\partial f(\x) := \{ \g\in \mathbb{R}^m| f(\y) - f(\x)  \ge  \g\T(\y - \x), \forall \y\in  \Rm \}$ denotes the subdifferential of $f$ at $\x$. 
	%
	The projection of $\x\in  \mathbb{R}^m$ on a nonempty closed convex set  $X\subseteq \mathbb{R}^m$ is denoted by $P_X(\x)$, i.e., $ P_X(\x) = \arg \min_{\y \in X} \|\x-\y\|$. 
	We denote by $\dist(\x,X)$ the (Euclidean) distance from $\x$ to  $X$, i.e., $\dist(\x,X) = \| \x -  P_X(\x)\|$. 
	The following inequality is called the nonexpansiveness property (see, e.g., \cite{Bertsekas99Book}):
	\begin{equation} \label{eqNonExpansive}
	\| P_X(\x) - P_X(\y)\| \le \| \x-\y\|, \quad \forall \x,\y \in \mathbb{R}^m. 
	\end{equation}
	We will employ the notion of regularity of the constraint sets, which plays an important role in the study of projection algorithms. This notion involves upper estimating the distance from a point to the intersection of a collection of closed convex sets in terms of the distance to each set (see \cite{Hoffmann92distance, Bauschke96projection}). Recalled next is the definition needed here, stated for a finite dimensional setting.
	\begin{defn}\label{defiRegularityConstraints}
		A collection of closed convex sets $\{X_i, i\in \mathcal{V}\}$ (with nonempty intersection) is \emph{regular} with respect to a set $B\subseteq \mathbb{R}^m$ if 
		$\exists~r_B \ge 1$ 
		such that $\forall \x \in B$
		\begin{equation}
		\textstyle \dist(\x, \cap_{i\in \mathcal{V}}X_i) \le r_B \max_{i\in \mathcal{V}} \dist(\x, X_i).	
		\label{eqRegularity}
		\end{equation}
	\end{defn}
	%
	\section{Problem statement and proposed algorithms} \label{secProblem}
	
	Consider a network of $n$ agents that can exchange information over a directed graph $\mathcal{G} = (\mathcal{V}, \mathcal{E})$. 
	All the agents share the objective of solving
	\begin{equation}
	\begin{split}
	\textstyle \min_{\x \in X} \quad &F(\x) := \textstyle\sum_{i\in \mathcal{V}} f_i(\x) \\
	\text{s.t.} \quad & X = \textstyle\bigcap_{i\in \mathcal{V}} X_i
	\end{split}
	\label{eqProblem}
	\end{equation}
	where each $f_i: \mathbb{R}^m \to \mathbb{R}$ is a convex function that we view as the private objective of agent $i$, and each $X_i \subseteq \mathbb{R}^m$ is a convex constraint set only known to agent $i$. Obviously, $F$ is also convex. 
	Let $F^*$ and $X^*$ denote the optimal value and the optimal solution set of \eqref{eqProblem}. 
	Let $U$ denote the convex hull of $\cup_{i\in \mathcal{V}} X_i$, i.e.,
	\begin{equation}
	U:=\textstyle \co\big(\bigcup_{i\in \mathcal{V}}X_i \big)
	\label{Defining_U}
	\end{equation}
	The following assumptions are adopted in the sequel. 
	
	\begin{assumption}\label{asmp_X}
		Problem \eqref{eqProblem} satisfies the following:\\[5pt]
		(a)~~\emph{(Constraint sets)} $X_i\subseteq \mathbb{R}^m$ are closed and convex and have a nonempty intersection (i.e., $X\neq \varnothing$), and $\{X_i, i\in \mathcal{V}\}$ is regular with respect to $U$.\\[5pt]
		(b)~~\emph{(Bounded subgradients)} For any $i\in \mathcal{V}$, $f_i$ is convex with subgradients bounded on $U$, i.e., 
		\begin{equation}
		\exists L\in (0,\infty), ~\|\g_{i}\| \le L, ~ \forall \g_i \in \partial f_i(\x),~ \forall \x\in U
		\label{eqBoundGi}
		\end{equation} 
		(c)~~The solution set $X^*$ is nonempty.
	\end{assumption}
	
	Here, the regularity assumption on $X_i$ is milder than imposing boundedness, allowing us to consider a broader class of optimization problems. 
	This assumption holds trivially  when the $X_i$ are identical; in fact, $r_U=1$ is a valid regularity constant (as per Definition~\ref{defiRegularityConstraints}). An unconstrained problem is a special case with $X_i = \Rm, \forall i\in \mathcal{V}$. 
	The regularity assumption is also satisfied if $X_i$ are compact and $X$ has a nonempty interior, i.e., $\exists \x_0 \in X$ and $\delta >0$ such that $\{\x: \|\x-\x_0\|\le \delta \} \subset X$; such assumptions are used in \cite{Nedic10AC,Lin16distributed}. In fact, it follows from  \cite[Cor.~2]{Hoffmann92distance} that $r_U = \sum_{i=1}^n (D_U/\delta)^i$ is a regularity constant, where $D_U$ denotes the diameter of $U$. 
	Other important cases include when $X_i$ are hyperplanes or half-spaces (see \cite{Bauschke96projection}). 
	Note that \cite{Lee16AC} assumes that there exists $c>0$ such that $\dist(\x, X) \le c~\dist(\x, X_i)$ for  $\forall i\in \mathcal{V}, \forall \x \in \mathbb{R}^m$, which is clearly stronger than our regularity condition. 

	%
	%
	Note also that 
	Assumption \ref{asmp_X}(b) 
	implies that each $f_i$, possibly non-differentiable, is $L$-Lipschitz continuous on $U$. 
	This assumption is not uncommon in the literature on non-smooth optimization and is satisfied by many cost functions; see, e.g.,  \cite{Nesterov04, Nedic10AC,Lin16distributed}.
	E.g., it holds for any convex functions if  $X_i$ are compact, since then $U$ is also compact. 
	When some $X_i$ are unbounded, this assumption includes, e.g., affine functions, the $\ell_1$ norm in the basis pursuit problem in compressed sensing, and $\ell_1$ regularized logistic regression in machine learning, and excludes certain classes of cost functions $f_i$ such as quadratic costs. 
	Note in passing that the recent works  \cite{LinRen17arXiv,Lee16AC} consider unbounded constraint sets but \cite{LinRen17arXiv}  assumes that  each $f_i$ is differentiable and has a nonempty and bounded set of critical points while  \cite{Lee16AC} considers differentiable functions with bounded and Lipschitz continuous gradients. 
	Finally, Assumption \ref{asmp_X}(c) can be satisfied when, e.g., at least one of the $X_i$ is compact, since then $X$ is compact. 
	In general, however, we do not require compactness of the constraint sets. 
	
	In our setting, agent $i$ only has access to $f_i$ and local information on its neighbors' broadcast state values, and a central coordinating node is absent. 
	Thus, the agents must collaborate in a distributed manner to solve problem \eqref{eqProblem}. 
	This involves local iterative computation along with  information diffusion. 
	We are interested in the scenario where the communication graph $\mathcal{G}$ 
	is directed and fixed. 
	We make the following additional blanket assumptions. 
	\begin{assumption} \label{asmpNetwork} 
		$\mathcal{G}$ is fixed and strongly connected. 
	\end{assumption}
	\begin{assumption}  \label{asmpIdentity} 
		The agents are labeled $1,2,..., n$ and their messages carry a unique identifier of the sender. Further, the agents know the value $n$ (or an upper bound).
	\end{assumption}

	Note that the assumption on unique identifiers is absent in most of the previous literature; it can be satisfied when the agents have different addresses. This is 
	usually the case in computer networking and telecommunication (e.g., when medium-access control addresses are used for packet collision avoidance), 
	especially in Ethernet and wireless networks.
	
	At any time instant, each agent exchanges its current state with its neighbors (in accordance with the directed network structure).  
	Upon receiving the information from its neighbors (including itself), agent $i$ updates its state in a weighted averaging scheme. 
	Thus, each edge $(i,j) \in \mathcal{E}$ is associated with a weight $w_{ij} \ge 0$ (locally chosen by agent $i$). 
	Let the weight matrix $W = [w_{ij}]$ satisfy the following condition. 
	
	\begin{assumption} \label{asmpWeight} 
		$W$ is row stochastic and satisfies $w_{ii}\!>\!0$ for $\forall i\!\in\! \mathcal{V}$, $w_{ij} \!>\! 0$ for $(i,j) \!\in\! \mathcal{E}$ and $w_{ij} \!=\! 0$ otherwise. 
	\end{assumption}
	
	
	Again we stress that unlike the case with existing algorithms in the literature, $W$ is only assumed to be row stochastic, and not either doubly stochastic or column stochastic. 
	Thus, each agent $i$ controls the $i$-th row of $W$ independently from the other rows,  
	giving it freedom in assigning the weights on its neighbors' information. 
	This explains why row stochastic matrices are more suitable for {\it ad hoc} wireless networks.%

	
	We now introduce the following distributed algorithm to solve problem \eqref{eqProblem} under all the assumptions above. 
	
	\begin{alg} \label{algSubgrad1}
		At $t=0$, agent $i$ initializes $\x_{i}(0)$ and sets $\z_{i}(0) \!=\! \e_i \!\in\! \mathbb{R}^n$ (or $\in  \mathbb{R}^{\tilde{n}}$ if only a bound $\tilde{n}$ on $n$ is known). For each $t \ge 0$, all agents update their states as follows:
		\begin{align} \label{eq_x}
		\x_{i}(t+1) &= \textstyle P_{X_i}\big( \sum_{j \in \mathcal{V}} w_{ij}\x_{j}(t) - \alpha(t)\frac{\g_{i}(t)}{z_{ii}(t)} \big) \\
		\z_{i}(t+1) &= \textstyle \sum_{j \in \mathcal{V}} w_{ij}\z_{j}(t),
		\label{eq_z}
		\end{align}
		where  
		$\g_{i}(t) \!\in\! \partial f_i(\sum_{j \in \mathcal{V}} w_{ij}\x_{j}(t))$,
		$\z_{i}(t) = [z_{i1},..., z_{in}]\T$ for each $i\in \mathcal{V}$, and $\alpha(t)$ is a step size (discussed later). 
	\end{alg}

	
	
	In essence,  \eqref{eq_x} is a modified version of the distributed projected subgradient (DPS) method \cite{Nedic10AC} where each private cost function's subgradient is scaled by $z_{ii}(t)$, 
	which is updated through \eqref{eq_z}. 
	Here,  \eqref{eq_z} is, in fact, a consensus iteration aiming to provide each agent $i\in \mathcal{V}$ with an estimate of  $\boldsymbol\pi  = [\pi_1, \ldots, \pi_n]\T$- the left normalized Perron eigenvector of $W$, i.e., the left eigenvector $\piB$ satisfying $\1\T \boldsymbol\pi  = 1$. This iteration resembles those used in \cite{QLL11, PGMS13}. 
	Of course, if each agent $i\in \mathcal{V}$ knows the $\pi_i$ in advance, then iteration \eqref{eq_z} is unnecessary as all the agents can simply use $z_{ii}(t) = \pi_i, \forall t\ge 0$. 
	(In fact, if initialized with $\z_{i}(0) = \boldsymbol{\pi}$, then it follows from \eqref{eq_z} that $\z_{i}(t) = \boldsymbol{\pi}$ for all $t\ge 0$.)
	In this case, our rescaling subgradient technique reduces to the reweighting scheme used in \cite{Tsianos12ACC,Lin16distributed}.
	
	Note that the DPS method in \cite{Nedic10AC} can be applied to time-varying networks but requires the weight matrix to be doubly stochastic for all time. 
	Further, for nonidentical constraints $X_i$, \cite{Nedic10AC} only considers complete graphs and assumes that $X$ has nonempty interior. 
	Later, \cite{Lin16distributed} extended the method to directed time-varying graphs possibly with (fixed and uniform) communication delays but still required doubly stochastic weight matrices and compact $X_i$ with nonempty interior. 
	Thus, the results in these works are not readily applicable to cases where the $X_i$ are unbounded  and/or $X$ has an empty interior (e.g., an $X_i$ includes linear equality constraints) and $W$ is only row stochastic. 
	Another extension in \cite{Nedic15AC} dealing with the unconstrained case employs column stochastic matrices. 
	Algorithm \ref{algSubgrad1} can be seen as an extension of DPS under the fixed network setting where only row stochastic weight matrices are used. 
	Note also that we assume the network is fixed during a run (i.e., from start to finish) of the algorithm. In case the network has to perform multiple runs (e.g., with new data) during its lifetime, we allow the network structure to change between any two consecutive runs, and our algorithm need not be adjusted except each agent $i$ may need to reselect new weights $w_{ij}$ for its (possibly) new neighbor set--a trivial task. 

	We also consider the following variation on Algorithm \ref{algSubgrad1}. 
	\begin{alg} \label{algSubgrad2}
		With the same initializations as in Algorithm \ref{algSubgrad1}, all agents use $\g_{i}(t) \!\in\! \partial f_i(\x_{i}(t))$ and perform
		\begin{align} \label{eq_x2}
		\x_{i}(t+1) &= \textstyle P_{X_i}  \Big( \sum_{j \in \mathcal{V}} w_{ij} \big( \x_{j}(t) - \alpha(t)\frac{\g_{j}(t)}{z_{jj}(t)} \big) \Big)  \\
		\z_{i}(t+1) &= \textstyle \sum_{j \in \mathcal{V}}w_{ij}\z_{j}(t),
		\label{eq_z2}
		\end{align}
	\end{alg}
	Note that \cite{Johansson08CDC} shows that the optimization and consensus steps in the usual DPS method can be interchanged, which, if a constant step size is used, often gives better convergence speed to a solution neighborhood \cite{Matei12}. 
	In-depth comparison between Algorithms~\ref{algSubgrad1} and \ref{algSubgrad2}, however, is out of the scope of this paper.%
	
	
	In this work, diminishing step size sequences satisfying the following assumption will be used to ensure convergence of our algorithms to the optimal solution. Convergence rate analysis will be performed under a less restrictive assumption. 
	\begin{assumption}\label{asmpStepSize}
		The sequence $\{\alpha(t)\}$ is positive nonincreasing with $\sum_{t=0}^{\infty}\alpha(t) = \infty$ and $\sum_{t=0}^{\infty}\alpha^2(t) < \infty$. 
	\end{assumption}
	There are many ways to choose the step size sequence $\alpha(t)$ satisfying this assumption, e.g., $\alpha(t) = \frac{c}{t^\gamma}, \forall t\ge 1$, for constants $c>0$ and $\gamma \in (0.5,1]$.
	
	\section{Basic relations and convergence results} \label{secMainResults}
	In this section, we prove the convergence of Algorithms 1 and 2.  
	We begin with a few basic results.
	
	Besides nonexpansivity \eqref{eqNonExpansive}, another property of a projection operator is given in the following lemma. 
	\begin{lemma}\label{lemProjectionInequality}(\cite{Nedic10AC})
		If $Y \subseteq \Rm$ is closed and convex, then $\| P_Y(\x)-\y \|^2 \!\le\! \|\x-\y\|^2 - \| P_Y(\x)-\x \|^2$, $\forall \x \!\in\! \Rm, \y\!\in\! Y$.
	\end{lemma}
	
	The next lemma follows from convexity of $\|\cdot\|^2$. 
	\begin{lemma}\label{lemJensen}
		If $\{a_i\}_{i=1}^n \subset \mathbb{R}_{+}$ satisfies $\sum_{i=1}^{n}a_i=1$, then $\|\sum_{i=1}^{n}a_i\x_i\|^2 \le \sum_{i=1}^{n}a_i \|\x_i\|^2$ for $\forall \{\x_i\}_{i=1}^n \subset \Rm$. 	
	\end{lemma}
	
	Convergence  of the power iteration of the weight matrix is recorded next, a consequence of the 
	Perron-Frobenius theorem (see \cite[Sec.~8.5]{Horn85} for details). 
	
	\begin{lemma} \label{lemWeight} 
		Let Assumptions \ref{asmpNetwork} and \ref{asmpWeight} 
		hold. Then $\lim_{t\to\infty} W^t = \1 \boldsymbol\pi\T$, 
		where $\boldsymbol{\pi} >\0$ is the normalized left Perron eigenvector of $W$. Moreover, the convergence is geometric with rate $\lambda \in (|\lambda_2(W)|,1)$, where $\lambda_2(W)$ is the second largest eigenvalue of $W.$
	\end{lemma}
	
	The next result, on convergence of the estimation step in \eqref{eq_z}, follows directly from the foregoing lemma; this result will be used in the sequel. 
	\begin{proposition}\label{propRateBound} \emph{(Convergence of $z_{ii}$)}
		Let Assumptions \ref{asmpNetwork}--\ref{asmpWeight} hold. Then for each $\lambda \in (|\lambda_2(W)|,1)$, there exists $C = C(\lambda,W)>0$ such that
		\begin{equation}
		| [W^t]_{ji} - \pi_i| \le C\lambda^t, \quad |z_{ii}(t) - \pi_i| \le C\lambda^t \label{eqRateBound}
		\end{equation} 
		for $\forall i,j \in \mathcal{V}$ and $\forall t$. Further, there exists $\eta >0$ satisfying  
		\begin{equation}
		\eta^{-1} \le z_{ii}(t) \le 1,\quad \forall t\ge 0,~ \forall i\in \mathcal{V}. \label{eqBoundZi}
		\end{equation}
	\end{proposition} 
	{\bf Proof.~}
		Let $Z(t) = [\z_1(t), \z_2(t),\cdots, \z_n(t)]\T$. 
		It follows from Algorithm \ref{algSubgrad1} that for any $t\ge 0$, $Z(t+1) = W Z(t)$ with $Z(0) = I$. Thus, $Z(t) = W^t, \forall t\ge 0$. Hence, \eqref{eqRateBound} follows by Lemma \ref{lemWeight} for some $C>0$ and $\lambda \in  (|\lambda_2(W)|,1)$. 
		
		Next, from \eqref{eq_z}, we have $z_{ii}(t+1) = \sum_{j \in \mathcal{V}}w_{ij}z_{ji}(t)$, $\forall i\in \mathcal{V}$, where $z_{ii}(0)=1, z_{ji}(0)=0, \forall j\neq i$. Clearly,  $1 \ge z_{ij}(t) \ge 0, \forall i,j\in \mathcal{V}, \forall t\ge 0$. 
		Since $\lim_{t\to \infty}z_{ii}(t) = \pi_i > 0$, there exists $t_0 \ge 0$ such that $z_{ii}(t) \geq \pi_i/2, \forall i\in \mathcal{V}, \forall t > t_0$. 
		Moreover, we have that $z_{ii}(t_0) \ge w_{ii}z_{ii}(t_0-1) \ge \ldots \ge w_{ii}^{t_0}z_{ii}(0) >0$ since $w_{ii}>0$ (cf. Assumption \ref{asmpWeight}). 
		Therefore,  $z_{ii}(t) > 0$ for any $t\in [0,t_0]$. 
		By taking 
		\begin{align} 
		\eta^{-1} = \min\{z_{ii}(t), \pi_i/2, \forall i\in \mathcal{V}, \forall t\in [0,t_0]\},\label{eqEta}
		\end{align}
		\eqref{eqBoundZi} follows as desired. \hfill$\blacksquare$
	
	\begin{remark}
		In the sequel, the parameters $C$, $\lambda$ and $\eta$ refer to the constants in Proposition \ref{propRateBound}.
	\end{remark}
	
	We now turn to \eqref{eq_x} and \eqref{eq_x2}. 
	Our next result describes a general relation on the overall evolution of the agents' states in terms of their distances from any $\vv\in X$ and the weighted averaged state vector $\bar{\x}(t)$, defined as
	\begin{equation}
	\bar{\x}(t) := \textstyle\sum_{j\in \mathcal{V}} \pi_j\x_{j}(t), \quad \forall t\ge 0. \label{eqXbar}
	\end{equation}
	The relation also involves the step size sequence $\alpha(t)$ and an error term $\big(F(\bar{\x}(t)) - F(\vv)\big),$ which in general is not the global objective error since $\bar{\x}(t)$ may not be in $X$; it is so if the constraint sets $\{ X_i \}_{i \in \mathcal{V}}$ are identical. 

	\begin{theorem} \label{thm1} \emph{(Bound on evolution of $\x_i$)}
		Let Assumptions \ref{asmp_X}--\ref{asmpWeight}
		be satisfied. Then for each of Algorithms \ref{algSubgrad1} and \ref{algSubgrad2}, the following holds for any $\vv\in X$ and $t\ge 0$:
		\begin{align}
		&\textstyle \sum_{i=1}^n \pi_i \|\x_{i}(t+1)-\vv\|^2  \nnb\\
		&\textstyle \le  (1+D_1\lambda^{2t})\sum_{i=1}^n \pi_i \|\x_{i}(t)-\vv\|^2 \nnb\\
		&~~~ -2\alpha(t)\big(F(\bar{\x}(t)) - F(\vv)\big) - \textstyle \sum_{i=1}^n \pi_i \| \phi_i(t) \|^2 \nnb\\
		&~~~ +D_2\alpha(t)\textstyle \sum_{i=1}^n \pi_i  \|\x_{i}(t)-\bar{\x}(t)\| + D_3\alpha^2(t),
		\label{eqThm1}
		\end{align}
		where 
		$D_1 = nCL\eta, D_2 = 2L\eta$, $D_3= L^2\eta^2+nLC\eta$,  and 
		\begin{align}
		\!\! \phi_i(t) := 
		\x_i(t+1) - \textstyle \big( \sum_{j=1}^n w_{ij}\x_{j}(t) - \alpha(t)\frac{\g_{i}(t)}{z_{ii}(t)} \big) \label{eqProjErr1}
		\end{align} 
		for Algorithm \ref{algSubgrad1}, while for Algorithm \ref{algSubgrad2}, $\phi_i(t)$ is defined as
		\begin{align}
		\!\!\phi_i(t) \!:=\! 
		\x_i(t\!+\!1) \!-\! \textstyle \sum_{j=1}^n w_{ij} \big(  \x_{j}(t) - \alpha(t)\frac{\g_{j}(t)}{z_{jj}(t)} \big) \label{eqProjErr2}
		\end{align} 
	\end{theorem}
	
	\vspace{-3pt}
		{\bf Proof.~}
		We prove the statement only for Algorithm \ref{algSubgrad1}; it follows for Algorithm~\ref{algSubgrad2} by a similar argument. 
		Let $\y_i(t) \!:=\! \textstyle\sum_{j=1}^n w_{ij}\x_{j}(t)$ and  
		$\vv\in X$. Then 
		\begin{align}
		&\|\x_{i}(t+1)-\vv\|^2  = \textstyle \big\| \y_i(t)- \vv -\alpha(t)\frac{\g_{i}(t)}{z_{ii}(t)} + \phi_i(t) \big\|^2  \nnb\\
		&\stackrel{\text{Lem.} \ref{lemProjectionInequality}}{\le} \textstyle \big\| \y_i(t)- \vv -\alpha(t)\frac{\g_{i}(t)}{z_{ii}(t)}\big\|^2 - \|\phi_i(t)\|^2   \label{eqThm1_0}.
		\end{align} 
		The first term on the right side of \eqref{eqThm1_0} equals
		\begin{align}
		\textstyle \!\!\!\!\|\y_i(t) \!-\! \vv\|^2 \!+\! \frac{2\alpha(t)}{z_{ii}(t)}\g_{i} (t)\T \!  (\vv \!-\! \y_i(t))  
		\!+\! \frac{\alpha^2(t)}{z_{ii}^2(t)}\|\g_{i}(t)\|^2. 
		\label{eqThm1a}
		\end{align} 
		We now derive an upper bound for each term in \eqref{eqThm1a}. Since $\y_i(t) \!-\! \vv \!=\! \sum_{j \in \mathcal{V}} w_{ij}(\x_j(t) \!-\! \vv)$, it follows that
		\begin{align}
		\|\y_i(t) - \vv\|^2 \stackrel{\text{Lem}. \ref{lemJensen}}{\le} \textstyle\sum_{j \in \mathcal{V}} w_{ij} \|\x_{j}(t) - \vv\|^2.  
		\label{eq1term}
		\end{align}
		Next, ignoring 
		$\frac{2\alpha(t)}{z_{ii}(t)}$, the second term in \eqref{eqThm1a} satisfies
		\begin{align}
		&\g_{i}(t)\T(\vv - \y_{i}(t)) 
		\le f_i(\vv) - f_i(\y_{i}(t))\nnb\\
		&\le f_i(\vv)  - f_i(\bar{\x}(t)) + \big| f_i(\y_{i}(t)) - f_i(\bar{\x}(t)) \big|	\nnb\\
		&\le  \textstyle   f_i(\vv)  - f_i(\bar{\x}(t))  + \! L \! \sum_{j \in \mathcal{V}} \! w_{ij} \left\| \x_{j}(t) \!-\! \bar{\x}(t) \right\|, 
		\label{eq2term}
		\end{align}
		where the first inequality holds since $\g_{i}(t) \in \partial f_i(\y_{i}(t))$, the second follows from the triangle inequality, 
		and the last from 
		Assumption \ref{asmp_X}(b) and the triangle inequality. 
		By continuing \eqref{eqThm1_0} and using \eqref{eqThm1a}, \eqref{eq1term}, \eqref{eq2term} and the conditions that $\|\g_{i}(t)\|\le L$ and $z_{ii}^{-1}(t)\le \eta, \forall i\in \mathcal{V}, \forall t\ge 0$, we have 
		$\|\x_{i}(t{+}1)-\vv\|^2\le \sum_{j \in \mathcal{V}}\! w_{ij} \|\x_{j}(t) \!-\! \vv\|^2 
		\!-\! \|\phi_i(t)\|^2 \!\!+\! \frac{2\alpha(t)}{z_{ii}(t)} \big(  f_i(\vv)  \!-\! f_i(\bar{\x}(t))  \big) 
		+ 2L\frac{\alpha(t)}{z_{ii}(t)} \sum_{j \in \mathcal{V}} w_{ij}\left\| \x_{j}(t) - \bar{\x}(t) \right\| 
		+ \alpha^2(t)L^2\eta^2$. 
		Thus, 
		\begin{align}
		&\!\!\! \textstyle\sum_{i\in \mathcal{V}} \pi_{i}\|\x_{i}(t\!+\!1)-\vv\|^2 \nnb\\
		&\!\!\! \le \textstyle \sum_{i\in \mathcal{V}} \pi_{i}\sum_{j \in \mathcal{V}} w_{ij} \|\x_{j}(t)- \vv\|^2 \nnb\\
		&\!\!\!	+ \textstyle 2 \!\sum_{i\in \mathcal{V}} \frac{\pi_{i} \alpha(t)}{z_{ii}(t)} \big(  f_i(\vv)  - f_i(\bar{\x}(t))  \big) -  \sum_{i\in \mathcal{V}} \pi_{i}\|\phi_i(t)\|^2 \nnb\\
		&\!\!\! + \textstyle 2L\!\sum_{i,j\in \mathcal{V}}\! \frac{\pi_{i} \alpha(t)}{z_{ii}(t)}  w_{ij}\!\left\| \x_{j}(t) \!-\! \bar{\x}(t) \right\| 	+\alpha^2(t)L^2\eta^2. 
		\label{eqThm1c}
		\end{align}
		Now consider each term on the right side of \eqref{eqThm1c}. First, 
		\begin{align}
		\!\!\!\!\textstyle \sum_{i,j\in \mathcal{V}} \pi_{i}w_{ij} \|\x_{j}(t)- \vv\|^2 = \sum_{i\in \mathcal{V}} \pi_{i}\|\x_{i}(t)-\vv\|^2, 
		\label{eqThm1c1}
		\end{align}
		where we have used the fact that $\boldsymbol\pi\T W=\boldsymbol\pi\T$. Second, the term $\sum_{i\in \mathcal{V}} \frac{\pi_{i}}{z_{ii}(t)} \big( f_i(\vv)- f_i(\bar{\x}(t))  \big)$ equals
		\begin{align}
		&\textstyle \!\sum_{i\in \mathcal{V}}\! f_i(\vv) - f_i(\bar{\x}(t)) 
		+\! \sum_{i\in \mathcal{V}}\! \frac{ \pi_{i} {-} z_{ii}(t)}{z_{ii}(t)} \big( f_i(\vv) - f_i(\bar{\x}(t)) \big) \nnb\\
		&\textstyle \!\le F(\vv) - F(\bar{\x}(t))   
		+ \sum_{i\in \mathcal{V}} \frac{|z_{ii}(t)- \pi_{i}|}{z_{ii}(t)} \left| f_i(\bar{\x}(t))- f_i(\vv) \right| \nnb\\
		&\le F(\vv)  - F(\bar{\x}(t))    +  nCL\eta \lambda^t \left\| \bar{\x}(t)- \vv \right\|, 
		\label{eqThm1c2}
		\end{align}
		where the last inequality follows from  $L$-Lipschitz continuity of $f_i$ and Proposition~\ref{propRateBound}. 
		Next, by using \eqref{eqBoundZi} and the relation $\boldsymbol\pi\T W=\boldsymbol\pi\T$ again, we have 
		\begin{align}
		\displaystyle \!\!\sum_{i,j\in \mathcal{V}}\! \frac{\pi_{i}w_{ij} }{z_{ii}(t)} \left\| \x_{j}(t) \!-\! \bar{\x}(t) \right\| 
		\le   \sum_{i\in \mathcal{V}} \eta \pi_{i} \!\left\| \x_{i}(t) \!-\! \bar{\x}(t) \right\|.
		\label{eqThm1c3}
		\end{align}
		Now, combining \eqref{eqThm1c}--\eqref{eqThm1c3} yields 
		\begin{align}
		&\!\!\!\!\textstyle  \sum_{i\in \mathcal{V}} \pi_{i}\|\x_{i}(t+1)-\vv\|^2  \nnb\\
		&\!\!\!\!\le \textstyle  \sum_{i\in \mathcal{V}} \pi_{i}\|\x_{i}(t)-\vv\|^2 
		- 2\alpha(t) \big( F(\bar{\x}(t))- F(\vv) \big)   \nnb\\
		& - \textstyle  \sum_{i\in \mathcal{V}} \pi_{i}\|\phi_i(t)\|^2  +2\alpha(t)nCL\eta \lambda^t \| \bar{\x}(t)-\vv \| \nnb\\
		& + \textstyle  2\alpha(t)L\eta \sum_{i \in \mathcal{V}} \pi_i \left\| \x_{i}(t) - \bar{\x}(t) \right\|  +\alpha^2(t)L^2\eta^2. 
		\label{eqThm1d}
		\end{align}
		Finally, by writing $ \bar{\x}(t)-\vv = \sum_{i\in \mathcal{V}} \pi_i (\x_i(t)-\vv) $ and then using ``$2ab\le a^2+b^2$"  and Lemma~\ref{lemJensen}, we have $2\alpha(t)\lambda^t \| \bar{\x}(t) \!-\! \vv \| \le \textstyle \alpha^2(t) \!+\! \lambda^{2t}  \sum_{i\in \mathcal{V}}\! \pi_i \| \x_i(t) \!-\! \vv \|^2$. 
		Using this bound for \eqref{eqThm1d} and then rearranging terms yields \eqref{eqThm1} as desired. \hfill$\blacksquare$
	
	It is worth highlighting the differences between this result, in particular \eqref{eqThm1}, with that obtained from the usual DPS method \cite{Nedic10AC} in the context of Algorithm \ref{algSubgrad1}. 
	First, since $\piB$ is nonuniform, we opt for employing the weighted average vectors $\bar{\x}(t)$ (and $\sum_{i\in \mathcal{V}} \pi_i \|\x_{i}(t)\!-\!\vv\|^2$) instead of the averages. 
	Second, the term $D_1\lambda^{2t}\sum_{i\in \mathcal{V}} \pi_i \|\x_{i}(t)\!-\!\vv\|^2$ (or more precisely  $2\alpha(t)nCL\eta \lambda^t \| \bar{\x}(t)-\vv \|$ in \eqref{eqThm1d}) arises as a result of each agent $i$ using an estimate $z_{ii}(t)$ of $\pi_i$ generated from the estimation step \eqref{eq_z}. 
	Finally, since $X_i$ may not be bounded or identical (or have a nonempty interior), the projection error $\phi_i$ is not guaranteed to be bounded a priori and $\big(F(\bar{\x}(t)) - F(\vv)\big)$ does not reflect a global objective error (as $\bar{\x}(t)$ need not be in $X$). 
	Therefore, quantifying the behaviors of these terms and errors will be the main challenging task in analyzing the convergence as well as the convergence rates of our  algorithms; this calls for new results that are more accessible than \eqref{eqThm1} which we develop in the sequel. 

	We now provide some  bounds on the terms  $\|\x_{i}(t){-}\bar{\x}(t)\|$ and $\|\phi_i(t)\|$ in \eqref{eqThm1} in terms of the step size sequence $\alpha(t)$ and the total projection error $\beta(t)$, defined as
	$$\textstyle\beta(t):=\sum_{i\in \mathcal{V}} \| \phi_i(t)\|, \quad\forall t\ge 0$$
	which, by the Cauchy-Schwarz inequality, satisfies  \begin{equation}
	\beta^2(t) \le n \textstyle\sum_{i\in \mathcal{V}} \|\phi_i(t)\|^2. \label{eqBoundBetaSquared}
	\end{equation}
	\begin{theorem}\label{thmXi_Xbar}
		Let Assumptions \ref{asmp_X}--\ref{asmpWeight} 
		hold. 
		Then  for both Algorithms \ref{algSubgrad1}~and~\ref{algSubgrad2}
		
		(a)~~ 	Let $D_4 := C\sum_{j \in \mathcal{V}} \|\x_{j}(0) \|$. For any $i\in \mathcal{V}$,
		\vspace{-2mm}
		\begin{align}
		\!\!\!\!\|\x_{i}(t)\!-\!\bar{\x}(t)\| 
		\!\le\! D_4 \lambda^t \!\!+  \!\!\!\!\!\!\sum_{0\le s\le t-1}\!\!\!\!\!\!\! \lambda^{t-1-s} \big( D_1\alpha(s) \!+\! C\beta(s) \big)
		\label{eqXi_Xbar_a}
		\vspace{-5mm}
		\end{align}
		%
		(b)~~ 	Let $\gamma(t) \!:=\! \alpha(t)\sum_{0\le s\le t-1}\lambda^{t-1-s} \beta(s)$ with $\gamma(0)\!=\!0.$
		If $\{\alpha(t)\}$ is nonincreasing, then 
		\begin{align}
		\gamma(t+1) \le \lambda\gamma(t) +  {\alpha(t)\beta(t)}. \label{eqBoundGamma_t}
		\end{align}
	\end{theorem}\vspace{-5mm}
	
	{\bf Proof.~}
		We only prove part (a); part (b) is straightforward and skipped for brevity. 
		First, we write \eqref{eq_x} and \eqref{eq_x2} 
		as $\x_{i}(t\!+\!1) \!=\! \textstyle\sum_{j \in \mathcal{V}} w_{ij}\x_{j}(t) + \boldsymbol\epsilon_{i}(t),$ 
		where $\boldsymbol\epsilon_i(t) \in \mathbb{R}^m$ is an error term. 
		As a result, $\x_{i}(t) = \textstyle\sum_{j \in \mathcal{V}} [W^t]_{ij}\x_{j}(0) + \sum_{0\le s\le t-1} \sum_{j \in \mathcal{V}} [W^{t-1-s}]_{ij}\boldsymbol\epsilon_{j}(s)$.  
		Since $\bar{\x}(t) \!=  \sum_{j \in \mathcal{V}} \pi_j\x_{j}(t) $ and $\boldsymbol\pi\T W \!= \boldsymbol\pi\T$, we have 
		$\bar{\x}(t) 
		=  \textstyle\sum_{j \in \mathcal{V}} \pi_j\x_{j}(0) + \sum_{0\le s\le t-1} \sum_{j \in \mathcal{V}} \pi_{j}\boldsymbol\epsilon_{j}(s)$. 
		As a result, $\|\x_{i}(t) - \bar{\x}(t)\| = \big\| 
		\sum_{j \in \mathcal{V}} ( [W^t]_{ij} \!-\! \pi_j ) \x_{j}(0) 
		\!+\! \sum_{s=0}^{t-1} \sum_{j \in \mathcal{V}} \!\! \left( [W^{t-1-s}]_{ij} \!-\! \pi_j \right)\! \boldsymbol\epsilon_{j}(s) \big\| 
		\le \sum_{j \in \mathcal{V}} | [W^t]_{ij} \!-\! \pi_j| \! \|\x_{j}(0) \| 
		+  \sum_{0\le s \le t-1} \sum_{j \in \mathcal{V}} \! | [W^{t-1-s}]_{ij} \!-\! \pi_j |\! \|\boldsymbol\epsilon_{j}(s)\|$. 
		Hence, by using the bound in \eqref{eqRateBound}, we have
		\begin{align} 
		\!\!\! \|\x_{i}(t) \!-\! \bar{\x}(t)\| \!\le\! 
		D_4 \lambda^t 
		\!+\! C\!\!\!\!\sum_{0\le s\le t-1}\!\!\!\lambda^{t-1-s}  \sum_{j \in \mathcal{V}} \|\boldsymbol\epsilon_{j}(s)\|.
		\label{eqProofLem_conv}
		\end{align}
		Now consider Algorithm \ref{algSubgrad1}, where  it follows from \eqref{eq_x} and \eqref{eqProjErr1} that 
		$ \boldsymbol\epsilon_i (t)= \phi_i(t) - \alpha(t)\frac{\g_{i}(t)}{z_{ii}(t)}.$	By using the triangle inequality and the facts that $\|\g_i(t)\| \le L$ (cf. Assumption \ref{asmp_X}(b)) and that $z_{ii}^{-1}\le \eta$ (see \eqref{eqBoundZi}), we obtain
		\begin{align}
		\| \boldsymbol\epsilon_i(t)\| \le \|\phi_i(t)\| + \alpha(t)L\eta, \quad \forall i\in \mathcal{V}. \label{eqErrorBound1}
		\end{align}
		Next, we show that this bound also holds for Algorithm \ref{algSubgrad2}. From \eqref{eq_x2} and \eqref{eqProjErr2} we have $ \boldsymbol\epsilon_i (t)= \phi_i(t) - \alpha(t)\sum_{j\in \mathcal{V}}w_{ij} \frac{\g_{j}(t)}{z_{jj}(t)}.$ As a result, for $\forall i\in \mathcal{V}$, 
		$\| \boldsymbol\epsilon_i(t)\| \!\le\! \|\phi_i(t)\| \!+\! \alpha(t)\!\sum_{j\in \mathcal{V}}\!w_{ij} \frac{\|\g_{j}(t)\|}{|z_{jj}(t)|}
		\!\le\! \|\phi_i(t)\| \!+\! \alpha(t)L\eta.$ 
		Finally, \eqref{eqXi_Xbar_a} follows from \eqref{eqErrorBound1} and \eqref{eqProofLem_conv}.%
		\hfill$\blacksquare$
		

	We note the following. First, \eqref{eqXi_Xbar_a} shows that the effect of initial conditions on the differences between agents' states vanishes exponentially. Second, we can view the last term on the right side of \eqref{eqXi_Xbar_a} as the sum of the convolutions of $\alpha(t)$ and $\beta(t)$ with $\lambda^t$. Thus, for convergence of the algorithms, we expect these terms to decay to zero under a suitable choice of $\alpha(t)$. E.g., when $\lim_{t\to \infty}\alpha(t)=0$, we show next that $\lim_{t\to\infty}\sum_{s=0}^{t-1}\! \lambda^{t-1-s}\alpha(s) = 0$. However, whether this also implies $\lim_{t\to\infty}\sum_{s=0}^{t-1}\! \lambda^{t-1-s}\beta(s) = 0$ is inconclusive since $\beta(t)$ depends on the agents' states and $X_i$. Finally, we introduced $\gamma(t)$  to study the behavior of the term $\alpha(t)\sum_{i\in \mathcal{V}} \pi_i  \|\x_{i}(t)-\bar{\x}(t)\|$ in \eqref{eqThm1}.

	\begin{corollary}\label{corLimBeta}
		In Theorem \ref{thmXi_Xbar}, if $\lim_{t\to \infty} \beta(t) = 0$, then $\lim_{t\to \infty} \gamma(t)=0$. Additionally, if  $\lim_{t\to \infty} \alpha(t) = 0$, then $\lim_{t\to \infty} \sum_{i\in \mathcal{V}} \pi_i  \|\x_{i}(t)-\bar{\x}(t)\|=0.$
	\end{corollary}
		{\bf Proof.~}
		Straightforward application of \cite[Lem.~7]{Nedic10AC}.\hfill$\blacksquare$

	Our next result is a consequence of  Theorems \ref{thm1} and \ref{thmXi_Xbar} under regularity of $X_i$.  
	Specifically, we apply the bounds in \eqref{eqXi_Xbar_a} and \eqref{eqBoundGamma_t} to \eqref{eqThm1}, and select coefficients to yield a more accessible relation, which will be key to proving convergence of the algorithms. 

	\begin{theorem}\label{thmMainBound2}
		Let Assumptions \ref{asmp_X}--\ref{asmpWeight} 
		hold. 
		If $\{\alpha(t)\}$ is nonincreasing, then for both 
		Algorithms \ref{algSubgrad1} and \ref{algSubgrad2}, we have
		\begin{align}
		&\textstyle\sum_{i\in \mathcal{V}} \pi_i \|\x_{i}(t+1)-\vv\|^2  + ab\gamma(t+1)\nnb\\
		&\le  \textstyle (1+D_1\lambda^{2t}) \sum_{i\in \mathcal{V}} \pi_i \|\x_{i}(t)-\vv\|^2 + ab\gamma(t) \nnb\\
		&\textstyle  -2\alpha(t)\big(F(\sB(t)) \!-\! F(\vv)\big) \!-\! D_6\! \sum_{i\in \mathcal{V}} \| \phi_i(t) \|^2 \!+\! D_{24} \alpha(t)\lambda^t \nnb\\
		&\textstyle   + D_{21} \alpha(t)\!\sum_{0\le s\le t-1}\! \lambda^{t-1-s}\alpha(s) + D'_3\alpha^2(t), 
		\label{eqMainBound2}
		\end{align}
		with $\sB(t) \!=\! P_X\big(\bar\x(t) \big)$, $\pi_{\min} \!=\! \min_{i\in \mathcal{V}}\pi_i$, $
		b \!=\! \sqrt{\frac{\pi_{\min}}{n}}, a =\frac{D_2'C}{(1-\lambda)b}, 
		D_2' \!=\! D_2 \!+\! \frac{2LR}{\pi_{\min}},$ $R$ a regularity constant of $\{X_i\}$,  
		$D_6 \!=\! \frac{\pi_{\min}}{2}, D_{24} \!=\! D_2'D_4, D_{21} \!=\! D_2'D_1$ and $ D'_3 \!=\! D_3+\! \frac{a^2}{2}$. 
	\end{theorem}
	%
		{\bf Proof.~}
		By adding and subtracting $F(\sB(t))$ and using Lipschitz continuity of $F$ we have 
		$F(\vv) - F(\bar{\x}(t)) \le F(\vv) - F(\sB(t)) + L \|\sB(t)- \bar{\x}(t)\|$. 
		Now we find an upper bound on the term $\|\sB(t)- \bar{\x}(t)\|$. 
		Since $\{X_i\}_{ i\in \mathcal{V} }$ is regular with respect to $U=\co(\cup_{i\in \mathcal{V}}X_i)$, it follows that 
		$\dist(\x, X) \le R \max_{i\in \mathcal{V}} \dist(\x, X_i)$, $\forall \x \in U$.
		Thus, 
		\begin{align}
		&\!\!\|\sB(t)- \bar{\x}(t)\| = \dist(\bar{\x}(t),X) 
		\le R \max_{i\in \mathcal{V}}\dist(\bar{\x}(t),X_i) \nnb\\
		&\!\!\!\le\! \sum_{i\in \mathcal{V}}\! \frac{R\pi_i}{\pi_{\min}} \dist(\bar{\x}(t),X_i) 
		\le \sum_{i\in \mathcal{V}}\! \frac{R\pi_i}{\pi_{\min}}\| \x_i(t) \!-\! \bar{\x}(t)\|, \label{eqBoundProjectionOnX}
		\end{align}
		where the last inequality holds since $\dist(\bar{\x}(t),X_i) \le \| \x_i(t) - \bar{\x}(t)\|$ (cf. Lem.~\ref{lemProjectionInequality}). 
		Hence, 
		$F(\vv) - F(\bar{\x}(t)) \le\! F(\vv) \!-\! F(\sB(t)) + \frac{L R}{\pi_{\min}}\sum_{i\in \mathcal{V}} {\pi_i} \|\x_i(t) -  \bar{\x}(t)\|. $
		Using this bound for \eqref{eqThm1}, we obtain 
		\begin{align}
		&\textstyle \sum_{i\in \mathcal{V}} \pi_i \|\x_{i}(t+1)-\vv\|^2  \nnb\\
		&\textstyle \le  (1+D_1\lambda^{2t})\sum_{i\in \mathcal{V}} \pi_i \|\x_{i}(t)-\vv\|^2 \nnb\\[-2pt]
		&\textstyle  -2\alpha(t)\big(F(\sB(t)) - F(\vv)\big) - \sum_{i\in \mathcal{V}} \pi_i \| \phi_i(t) \|^2 \nnb\\[-2pt]
		&\textstyle   +(D_2 \!+\! \frac{2LR}{\pi_{\min}})\alpha(t)\sum_{i\in \mathcal{V}} \pi_i  \|\x_{i}(t) \!-\! \bar{\x}(t)\| \!+\! D_3\alpha^2(t).  \nnb
		\end{align}
		Next, adding $ab\gamma(t+1)$ to both sides of this inequality and using bounds \eqref{eqXi_Xbar_a} and \eqref{eqBoundGamma_t}, we further have
		\begin{align}
		&\textstyle \!\!\!\!\!\sum_{i\in \mathcal{V}} \pi_i \|\x_{i}(t+1)-\vv\|^2  + ab\gamma(t+1)\nnb\\
		&\textstyle \!\!\!\!\!\le  (1+D_1\lambda^{2t})\sum_{i\in \mathcal{V}} \pi_i \|\x_{i}(t)-\vv\|^2 + ab\gamma(t) \nnb\\
		&\!\!\!\!\! + ab(\lambda-1)\gamma(t) + ab \alpha(t)\beta(t)\nnb\\
		&\textstyle \!\!\!\!\! -2\alpha(t)\big(F(\sB(t)) \!-\! F(\vv)\big) \!-\! \sum_{i\in \mathcal{V}} \pi_i \| \phi_i(t) \|^2 \!+\! D_3\alpha^2(t)\nnb\\
		&\!\!\!\!\! +D_{24}\alpha(t) \lambda^t \!+ D_{21} \alpha(t)\!\!\!\!\!\!\sum_{0\le s\le t-1}\!\!\!\!\!\! \lambda^{t-1-s}\alpha(s) \!+\! D_2'C\gamma(t). \label{eqMainbound2a}
		\end{align}
		Since $a =\frac{D_2'C}{(1-\lambda)b}$, the terms $ab(\lambda-1)\gamma(t)$ and $D_2'C\gamma(t)$ cancel out. Furthermore,  we also have 
		$		2ab\alpha(t)\beta(t) 
		\le a^2\alpha^2(t) + b^2\beta^2(t)
		\le a^2\alpha^2(t) + b^2n\sum_{i\in \mathcal{V}}\| \phi_i(t) \|^2,$
		where the last inequality follows from  \eqref{eqBoundBetaSquared}. 
		As a result, with $b^2=\frac{\pi_{\min}}{n}$, it can be verified that 
		$ab\alpha(t)\beta(t) \!-\! \sum_{i\in \mathcal{V}}\! \pi_i \| \phi_i(t) \|^2 
		\!\le\! \frac{a^2}{2}\alpha^2(t) \!-\!\frac{\pi_{\min}}{2}\!\sum_{i\in \mathcal{V}}\!\| \phi_i(t) \|^2.$ 
		It remains to apply the relations above to \eqref{eqMainbound2a} and then rearrange terms to obtain \eqref{eqMainBound2}. \hfill$\blacksquare$
	
	Note that \eqref{eqMainBound2} holds uniformly on $X$ since $D_i$ are independent of $\vv \!\in\! X$. When restricted to $X^*$, we immediately have a relation between 
	$\sum_{i\in \mathcal{V}} \pi_i \|\x_{i}(t)-\vv^*\|^2$ and the global objective error $F(\sB(t)) - F^*$, 
	both of which are desired to converge under a suitable choice of $\{\alpha(t)\}$.
	We are now ready to give a convergence result that applies to both Algorithms \eqref{eq_x}-\eqref{eq_z} and \eqref{eq_x2}-\eqref{eq_z2}, whose proof will be based on Theorem \ref{thmMainBound2} and the following lemma. 
	
	\begin{lemma}(\cite{RS71})\label{lem_superMartingale} 
		Let $\{v_t\}$, $\{u_t\}$, $\{b_t\}$ and $\{c_t\}$ be nonnegative sequences 
		with $\sum_{t=0}^{\infty} b_t \!<\! \infty$, $\sum_{t=0}^{\infty} c_t \!<\! \infty$ and
		\begin{equation}
		v_{t+1} \le (1+b_t)v_t- u_t + c_t, \quad \forall t\ge 0.
		\end{equation}
		Then $\{v_t\}$ converges and $\sum_{t=0}^{\infty} u_t < \infty$. 
	\end{lemma} 
	
	\begin{theorem} \label{thmConvergence} \emph{(Convergence)}
		Under Assumptions 
		\ref{asmp_X}-\ref{asmpStepSize},  
		both Algorithms \ref{algSubgrad1} and \ref{algSubgrad2} yield convergence to an optimal solution, i.e., $\exists \x^*\in X^*: \lim_{t\to\infty} \x_{i}(t) = \x^*, \forall i\in \mathcal{V}.$
	\end{theorem}
		{\bf Proof.~} We proceed in two steps: (i) apply Lemma \ref{lem_superMartingale}  to \eqref{eqMainBound2}, and then (ii) prove convergence.  
		
		{\bf Step~(i): } Let $\x^{\dagger}$ be arbitrary in $X^*$ and define the nonnegative sequences $\{v_t\}$, $\{u_t\}$, $\{b_t\}$ and $\{c_t\}$ as follows: 
		\begin{align}
		v_t &:= \textstyle \sum_{i\in \mathcal{V}} \pi_i \|\x_{i}(t)-\x^{\dagger} \|^2 + ab\gamma(t), \qquad b_t := D_1\lambda^{2t}, \nnb\\
		u_t &:= \textstyle 2\alpha(t)(F({\sB}(t)) - F^*) + D_6 \sum_{i\in \mathcal{V}} \| \phi_i(t) \|^2 , \nnb\\[-1pt]
		c_t &:= D_{24} \alpha(t)\lambda^t \!+ D_{21} \alpha(t)\!\!\!\!\sum_{0\le s\le t-1}\!\!\!\!\! \lambda^{t-1-s}\alpha(s) + D'_3\alpha^2(t).\nnb 
		\end{align}
		Adding the nonnegative term $D_1\lambda^{2t} ab\gamma(t)$ to the right side of \eqref{eqMainBound2} yields  
		$v_{t+1} \le (1+b_t)v_t - u_t + c_t, \forall t\ge 0.$
		To apply Lemma \ref{lem_superMartingale}, we show that $\{b_t\}$ and $\{c_t\}$ are summable. 
		Since $\lambda \in (0,1)$, we have $\sum_{t=0}^{\infty} b_t = D_1/(1-\lambda^2)$. Next, consider each term in $c_t$. First, $\sum_{t\ge 0}\alpha^2(t)\!<\!\infty$ by Assumption~\ref{asmpStepSize}. Second, using the fact $2\alpha(t)\lambda^t \!\le\! \alpha^2(t)\!+\!\lambda^{2t}$ yields $\sum_{t\ge 0} 2\alpha(t)\lambda^t \!\le\! \sum_{t\ge 0} \alpha^2(t) \!+\! \sum_{t\ge 0} \lambda^{2t} \!<\! \infty$.
		Third, by monotonicity of  $\{\alpha(t)\}$ (cf. Assumption~\ref{asmpStepSize}) the second term in $c_t$ satisfies: 
		$
		\alpha(t)\!\sum_{s = 0}^{t-1}\! \lambda^{t-1-s}\alpha(s) 
		\le \sum_{s = 0}^{t-1}\! \lambda^{t-1-s}\alpha^2(s).$ 
		Thus, for $\forall N\ge 1$,
		\begin{align}
		&\textstyle \sum_{t=1}^N \alpha(t) \sum_{s=0}^{t-1} \lambda^{{t-1}-s}\alpha(s) \nnb\\
		&\textstyle \le \sum_{0\le s \le t \le N-1}\! \lambda^{t-s}\alpha^2(s) 
		\le  \sum_{s=0}^{N-1}\!\alpha^2(s)\! \sum_{t \ge  s}\! \lambda^{t-s}\nnb\\
		&\textstyle = \sum_{0\le s \le N-1} \frac{\alpha^2(s)}{1-\lambda} \le  \frac{1}{1-\lambda} {\sum_{s\ge 0}\alpha^2(s)}<\infty.\label{eqThm2f}
		\end{align}
		Thus $\{c_t\}$ is summable. Therefore, by Lemma \ref{lem_superMartingale}, 
		there exists $\delta \ge 0$ such that
		\begin{align}
		&\textstyle \!\!\!\!\!\lim_{t\to\infty} \sum_{i\in \mathcal{V}} \pi_i \|\x_{i}(t)-\x^\dagger\|^2 + ab\gamma(t)= \delta 
		\label{eqThm2d}\\
		&\textstyle \!\!\!\!\!\sum_{t\ge 0}\! \Big[ \alpha(t)\big[ F(\sB(t)) \!-\! F^*\big] \!+\! \frac{D_6}{2}\!\sum_{i\in \mathcal{V}}\! \| \phi_i(t) \|^2 \Big] < \infty \label{eqThm2e}
		\end{align}
		%
		{\bf Step~(ii):} 
		First, by \eqref{eqThm2e}, we have $\lim_{t\to \infty} \sum_{i\in \mathcal{V}} \| \phi_i(t) \|^2 =0.$ Thus, 
		$\lim_{t\to \infty} \beta(t)=0$, which by Corollary \ref{corLimBeta} yields $\lim_{t\to \infty} \gamma(t)=0$. It then follows from \eqref{eqThm2d} that 
		\begin{align}
		\textstyle \lim_{t\to\infty} \sum_{i\in \mathcal{V}} \pi_i \|\x_{i}(t)-\x^\dagger\|^2= \delta. \label{eqThm2d2}
		\end{align}
		As a result, for each $i\in \mathcal{V}$, $\{\x_i(t)\}_{t\ge 0}$ is a bounded sequence. Thus so are $\{\bar{\x}(t)\}_{t\ge 0}$ and $\{ \sB(t)\}_{t\ge 0}$. 

		Next, since  $\sum_{t\ge0}\alpha(t) = \infty$, it then follows from \eqref{eqThm2e} that $\liminf_{t\to\infty} F(\sB(t)) = F^*.$ 
		Thus, there exists a subsequence $\{\sB(t_k)\} \subseteq \{\sB(t)\}$  such that 
		\begin{equation}
		\textstyle \lim_{k\to\infty} F(\sB(t_k)) = F^*. 
		\label{eqThm2e2}
		\end{equation}
		Since $\{\sB(t_k)\}$ is also bounded, it has a convergent subsequence $\{\sB(t_l)\} \!\subseteq\! \{\sB(t_k)\}$, i.e., $\lim_{l\to\infty} \sB(t_l) = \x^*$ for some $\x^*\in X$ (since $X$ is closed). We next show that $\x^*\in X^*$. 
		By continuity of $F$, 
		we have $\lim_{l\to\infty} F(\sB(t_l)) \!=\! F(\x^*)$, 
		which in view of \eqref{eqThm2e2} implies that $F(\x^*) \!=\! F^*$. 
		By convexity of $F$, we conclude that $\x^* \!\in\! X^*$. 
		Since $\x^\dagger \!\in\! X^*$ was chosen arbitrarily, we can let $\x^\dagger \!=\! \x^*$. 
		
		It now remains to show that $\delta =0$, which  by \eqref{eqThm2d2} will then complete the proof. By the triangle and Cauchy-Schwarz inequalities, 
		it can be verified that 
		$\| \x_{i}(t) \!-\!  \x^* \|^2 \!\le 3 (\| \x_{i}(t) \!-\! \bar{\x}(t)\|^2   \!+\! \|  \bar{\x}(t) \!-\! \sB(t)\|^2 \!+\! \| \sB(t) \!-\! \x^* \|^2)$
		Next, since  
		$\|\bar{\x}(t) - \sB(t)\| \le\! \frac{R}{\pi_{\min}}\sum_{i\in \mathcal{V}}\! {\pi_i}\| \x_i \!-\! \bar{\x}(t)\|$ (see   \eqref{eqBoundProjectionOnX}), we have 
		$\|\sB(t)- \bar{\x}(t)\|^2 \le\! \frac{R^2}{\pi_{\min}^2}\sum_{i\in \mathcal{V}}\! \pi_i\| \x_i \!-\! \bar{\x}(t)\|^2$ by Lemma~\ref{lemJensen}. 
		As a result, 
		$	\!\frac{\| \x_{i}(t) \!-\!  \x^* \|^2}{3} 
		\le 
		\| \x_{i}(t) \!-\! \bar{\x}(t)\|^2 \!+\! \frac{R^2}{\pi_{\min}^2}\sum_{i\in \mathcal{V}} \pi_i\| \x_i \!-\! \bar{\x}(t)\|^2  + \| \sB(t) - \x^* \|^2.$
		Multiplying both sides by $\pi_i$ then summing over $i\in \mathcal{V}$ yields the following, with $R'=1\!+\!{R^2}/{\pi^2_{\min}}$: 
		\begin{align}
		\sum_{i\in \mathcal{V}}\! \frac{\pi_i}{3} \| \x_{i}(t) \!-\!  \x^* \|^2 
		\!\le\!  R'\!\sum_{i\in \mathcal{V}}\! \pi_i \| \x_{i}(t) \!-\! \bar{\x}(t)\|^2   \!\!+\! \| \sB(t) \!-\! \x^* \|^2 \nnb
		\end{align}
		Taking $\liminf_{t\to\infty}$ of both sides and using \eqref{eqThm2d2} yields:
		\begin{align}
		\frac{\delta}{3}&\le \liminf_{t\to \infty} \Big(R'\sum_{i\in \mathcal{V}}\pi_i \| \x_{i}(t) - \bar{\x}(t)\|^2    + \| \sB(t) - \x^* \|^2\Big) \nnb \\
		&= \liminf_{t\to \infty} \| \sB(t) - \x^* \|^2. \label{eqThm2f4}
		\end{align}
		Here we have used the superadditivity property of $\liminf$ and the fact that $\lim_{t\to\infty} \sum_{i\in \mathcal{V}} \pi_i \| \x_{i}(t) - \bar{\x}(t)\|^2 = 0$ since  $\lim_{t\to\infty}\beta(t)=0$ (see Corollary \ref{corLimBeta}). Since   $\{\sB(t_l)\} \to \x^*$, we have $\liminf_{t\to \infty} \| \sB(t) - \x^* \| = 0$, which in view of \eqref{eqThm2f4} implies that $\delta = 0$. \hfill$\blacksquare$
	\section{Rate of convergence} \label{secRoC}
	We now discuss the convergence rate of our algorithms, which evidently depends on the choice of $\alpha(t)$. 
	Since the estimation step \eqref{eq_z} converges exponentially, one should expect that the convergence rate of the objective error is equivalent to that of usual distributed subgradient methods for the case when $X_i$ are identical and/or compact. 
	We emphasize, however, that such assumptions are relaxed in our work, i.e., the sets $X_i$ can be nonidentical and unbounded. Moreover, the global constraint set $X$ can also have empty interior. Thus, the agents'  estimates $\x_i(t)$ as well as their weighted average $\bar{\x}(t)$ need not be in the set $X$ at any time $t$. 
	As a result, local analysis around an optimal solution does not readily apply. 
	
	In this work, to quantify the distance from the optimum, we first propose to use a combined error term  involving (i) the distance from a local estimate $\tilde{\x}_i(t)$ of each agent to some $\tilde{\sB}(t) \in X$, and (ii) a global objective error evaluated at  $\tilde{\sB}(t)$, namely  $F(\tilde{\sB}(t)) - F^*$. 
	Specifically, define 
	\begin{align}
	\!\!\!\!\tilde{\x}_i(t) \!:=\! \frac{\sum_{k=0}^{t}\alpha(k)\x_i(k)}{\sum_{k=0}^{t}\alpha(k)},~~ 
	\tilde{\sB}(t) \!:=\! \frac{\sum_{k=0}^{t}\alpha(k)\sB(k)}{\sum_{k=0}^{t}\alpha(k)}.
	\end{align}
	Here, 
	$\tilde{\x}_i(t)$ can be computed locally by agent $i$ but might not lie in $X$. In contrast, $\tilde{\sB}(t)$ always belongs to $X$ but is not directly  available to each agent. 
	Our next theorem asserts that both errors $\| \tilde{\x}_i(t) - \tilde{\sB}(t)  \|$ and $F(\tilde{\sB}(t)) - F^*$ decay as $O\big({\sum_{k=0}^{t} \alpha^2(k)}/{\sum_{k=0}^{t} \alpha(k)} \big)$. 
	Based on this result, we then show that the same rate holds true for each agent's estimate, namely $|F(\tilde{\x}_i(t)) - F^*|$. These results are established based on the following lemma.%
	
	\begin{lemma} \label{lemIneqProd}
		Let $D >0$ and $\lambda\in (0,1)$ and 
		define $g_{t}(D,\lambda) := \prod_{0\le k\le t} (1+D\lambda^k)$. Then $\{g_t(D,\lambda)\}_t$ is positive, increasing and convergent. Moreover, 
		\begin{align}
		\textstyle 1+\frac{D}{1-\lambda} \le g_{\infty}(D,\lambda) \le e^{\frac{D}{1-\lambda}}. \label{ineqProd}
		\end{align}
	\end{lemma}
		{\bf Proof.~}
		For any $T\ge 1$, we have 
		$1+ D\!\sum_{0\le t\le T}\! \lambda^t \le g_T(D,\lambda)\le e^{D\sum_{t=0}^{T} \lambda^t},
		$
		where the second inequality follows from the basic relation that $1+x \le e^{x}$, $\forall x\ge 0$. Letting $T\to \infty$ yields \eqref{ineqProd}. The rest is obvious.  \hfill$\blacksquare$

	\begin{theorem} \label{thmConvergenceRate} \emph{(Convergence bound)}
		Let Assumptions 
		\ref{asmp_X}--\ref{asmpWeight} 
		hold and 
		let $\{\alpha(t)  \}$ be nonnegative and nonincreasing. 
		For both Algorithms \ref{algSubgrad1} and \ref{algSubgrad2}, the following holds:
		\begin{equation}
		C_0\| \tilde{\x}_i(t) \!-\! \tilde{\sB}(t)  \| \!+\!  F(\tilde{\sB}(t)) \!-\! F^* 
		\!\le\! E(t), \quad \forall t\ge 0, 
		\label{eqFConvergence}
		\end{equation}
		where $E(t) := \big({C_1+C_2\sum_{k=0}^{t}\alpha^2(k)}\big)/{\sum_{k=0}^{t}\alpha(k)}$ and 
		\begin{align}
		C_0 &= \frac{D_6(1-\lambda)}{{2n(Rn+1)C} },  
		C_1 \!=  R_1\!+\!\frac{R_2D_{24}\alpha(0)}{1-\lambda} \!+\! \frac{D_6D_{4}\alpha(0)}{2n C} \nnb\\
		C_2 &= R_2\big(\frac{D_{21}}{1-\lambda}\! + D'_3\big) + \frac{D_6}{2n}\big(  \frac{D_{1}}{ C} + \frac{1}{4}\big)\nnb
		\end{align}
		with $R_2 \!:=\! \frac{1}{2}g_{\infty}(D_1,\lambda^2)$ 
		and $R_1 \!:=\!  R_2\! \sum_{i\in \mathcal{V}}\! \pi_i \|\x_{i}(0){-}\x^* \|^2$  for any $\x^*\in X^*$. Here, $R$ is a regularity constant of $\{X_i\}_{ i\in \mathcal{V} }$ with respect to $U$. Moreover, 		
		\begin{align}
		|F(\tilde{\x}_i(t)) - F^*| \le E(t)\big( 1+{nL}/{C_0} \big). \label{eqFConvergenceIdenticalConstraints}
		\end{align}%
		
	\end{theorem}
		{\bf Proof.~}
		We proceed through the following three steps: 
		(i) Use the bound \eqref{eqMainBound2} in Theorem \ref{thmMainBound2} to upper estimate the sum $\!\sum_{k=0}^t\! 2\alpha(k)(F(\sB(k)) \!-\! F^*) \!+\! D_6 \sum_{i\in \mathcal{V}} \| \phi_i(k) \|^2$ in terms of $\sum_{k=0}^{t}\alpha(k)$ and  $\sum_{k=0}^{t}\alpha^2(k)$; (ii) Relate the left side of \eqref{eqFConvergence} to this sum by using the convexity of $F$ and the bounds given in Theorem \ref{thmXi_Xbar}; (iii) Prove  \eqref{eqFConvergenceIdenticalConstraints} using Lipschitz continuity of $F$ and~\eqref{eqFConvergence}.
		
		{\bf Step (i):}  
		Let 
		$\{v_t\}$, $\{u_t\}$, $\{b_t\}$ and $\{c_t\}$ be defined as in Step (i) of the proof of Theorem~\ref{thmConvergence}. 
		Further, set $$\Phi_t := \textstyle\sum_{i\in \mathcal{V}} \| \phi_i(t) \|^2.$$
		By using Theorem~\ref{thmMainBound2} and adding the nonnegative term $b_tab\gamma(t)$ to the right  side of \eqref{eqMainBound2}, we have  
		that $v_{t+1} \!\le\! (1+b_t)v_t   - u_t + c_t,~\forall t\ge 0$,
		which then implies that 
		\begin{align}
		\!\!v_{t+1} 
		\le\! \!\!\prod_{0\le k\le t}\!\!\!(1+b_k)v_0 &+ \!\!\!\sum_{0\le k\le t-1}\!\! (c_k-u_k)\!\!\!\!\!\prod_{k+1\le s\le t}\!\!\!\!\! (1+b_s) \nnb\\
		&+ (c_t-u_t).\label{eqBoundVt}
		\end{align}
		By Lemma \ref{lemIneqProd}, the following holds for any $t,k \ge 0$:
		$$
		\textstyle 1 < \prod_{k\le s \le t}(1+b_s) <  g_{\infty}(D_1,\lambda^2)
		=: D_e.
		$$
		As a result, \eqref{eqBoundVt} yields 
		\begin{align}
		\textstyle v_{t+1} \le D_ev_0 + \sum_{0\le k\le t} D_e c_k- \sum_{0\le k\le t} u_k. 
		\end{align}
		By rearranging terms and using the fact that $v_{t+1} \ge 0$, we have (recalling the definitions of $u_t$, $R_1$ and $R_2$)
		\begin{align}
		\!\!\!\sum_{0\le k\le t}\!\!\! \big[ \alpha(k)(F(\sB(k)) \!-\! F^*) \!+\! \frac{D_6}{2} \Phi_k \big] \!\le\! R_1 \!+\! R_2\!\!\!\sum_{0\le k\le t}\!\!\! c_k
		\label{eqSumUt}
		\end{align}
		Next, we will derive an upper bound on the term $\sum_{k=0}^{t} c_k$ based on the following estimates: 
		\begin{align}
		&\textstyle \sum_{0\le k\le t}\alpha(k) \lambda^k \le \sum_{0\le k\le t}\alpha(0) \lambda^k \le \frac{\alpha(0)}{1-\lambda}, \label{SumAkLamk}\\
		&\textstyle  \sum_{k=1}^{t} \alpha(k)\! \sum_{s=0}^{k-1} \lambda^{{k-1}-s}\alpha(s)
		\le  \sum_{s=0}^{t-1} \frac{\alpha^2(s)}{1-\lambda},\label{SumAkLamAs}
		\end{align}
		where \eqref{SumAkLamAs} is obtained from \eqref{eqThm2f}. 
		Hence, 
		$\!\sum_{k=0}^{t}\! c_k 
		\!\le\! \frac{D_{24}\alpha(0)}{1-\lambda}+ (\frac{D_{21}}{1-\lambda} + D'_3)\!\sum_{k=0}^{t}\!\alpha^2(k).$ 
		Therefore, 
		\begin{align}
		\!\!\sum_{0\le k \le t}\!\!\!\! \big[  \alpha(k)(F(\sB(k)) \!-\! F^*) \!+\! \frac{D_6}{2} \Phi_k \big]
		\!\le\! M_1 \!+\! M_2\!\!\!\!\sum_{0\le k\le t}\!\!\!\!\alpha^2(k) \label{eqFbarBound}
		\end{align}
		with 
		$
		M_1 =  R_1+\frac{R_2D_{24}\alpha(0)}{1-\lambda}$ and  $M_2 =  R_2\big(\frac{D_{21}}{1-\lambda}\! + D'_3\big).
		$
		
		{\bf Step (ii): }
		Now we derive lower bounds on the left side of \eqref{eqFbarBound}. 
		By convexity of $F$, 
		\begin{align}
		F(\tilde{\sB}(t)) - F^* \le \sum_{0\le k \le t}\frac{\alpha(k) \bigl( F(\sB(k)) -F^* \bigr)}{\sum_{\tau=0}^{t}\alpha(\tau)} . 
		\label{eqFEst1}
		\end{align}
		Next, we relate the term $\| \tilde{\x}_i(t) - \tilde{\sB}(t)  \| $ to $\sum_{k=0}^{t}\Phi_k.$
		By the triangle inequality, it can be shown that 
		\begin{align}
		\| \tilde{\x}_i(t) - \tilde{\sB}(t) \| 
		\le \sum_{0\le k \le t}\frac{\alpha(k) \| \x_i(k) - \sB(k) \| }{\sum_{\tau=0}^{t}\alpha(\tau)}.
		\label{eqFEst2}
		\end{align}
		We now quantify the numerator of the right side of \eqref{eqFEst2}. First, recall from  \eqref{eqXi_Xbar_a} that
		\begin{align}
		\!\|\x_{i}(t)-\bar{\x}(t)\| 
		\le D_4 \lambda^t \!+  \!\!\!\!\sum_{0\le s\le t-1}\!\!\!\!\! \lambda^{{t-1}-s} \big( D_1\alpha(s) + C\beta(s) \big).\nnb
		\end{align}
		Second, using the regularity assumption on $\{X_i\}_{ i\in \mathcal{V} }$ and basic properties of the $\dist(\cdot)$ function yields
		\begin{align}
		&\|\sB(k) - \bar{\x}(k)\| = \dist(\bar{\x}(k),X) 
		\le R \max_{i\in \mathcal{V}}\dist(\bar{\x}(k),X_i) \nnb\\
		&\textstyle \le R\sum_{i\in \mathcal{V}}\dist(\bar{\x}(k),X_i) 
		\le R\sum_{i\in \mathcal{V}} \|\x_i(k) - \bar{\x}(k)\|.\nnb
		\end{align}
		By the triangle inequality and the two previous relations, it can be shown that 
		\begin{align}
		\frac{\| \x_i(k) \!-\! \sB(k)\|}{(Rn+1)C} \!\le \frac{D_4}{C} \lambda^k  +\!\!\!\!   \sum_{0\le s\le k-1} \!\!\!\!\!\! \lambda^{{k-1}-s}\big(\frac{D_1\alpha(s)}{C}+\beta(s) \big)\nnb
		\end{align}
		which yields (see the definition of $\gamma(t)$ in Theorem~\ref{thmXi_Xbar}-b)
		\begin{align}
		&\!\! \frac{1}{{(Rn+1)C} }\sum_{0\le\! k\le t}\alpha(k){\| \x_i(k) \!-\! \sB(k)\|}\nnb\\
		&\!\! \le\! \frac{D_4}{C}\!\sum_{k=0}^{t}\!\alpha(k)\lambda^k \!+\! \frac{D_1}{C}\!\sum_{k=1}^{t}\!\alpha(k)\!\sum_{s = 0}^{k-1}\!  \lambda^{k-1-s}\alpha(s) 
		\!+\!\sum_{k=0}^{t}\!\gamma(k)\nnb\\
		&\!\!\!\!\stackrel{\eqref{SumAkLamk},\eqref{SumAkLamAs}}{\le} \!
		\frac{D_{4}\alpha(0)}{(1-\lambda)C} \!+ \frac{D_{1}}{(1-\lambda)C}\!\sum_{s=0}^{t}\!\alpha^2(s) \!+\! \sum_{k=0}^{t}\!\gamma(k).\!
		\label{eqSumAkXi_Xbar}
		\end{align}
		The sum $\sum_{k=0}^{t}\!\gamma(k)$ can be bounded as follows. By \eqref{eqBoundGamma_t} and noting that $\gamma(0)=0$ and $\gamma(t)\ge 0, \forall t\ge 1$, we have 
		\begin{align}
		\textstyle \sum_{k=0}^{t} \gamma(k)
		\le  \lambda \sum_{k=0}^{t-1} \gamma(k) + \sum_{k=0}^{t-1}\alpha(k)\beta(k).
		\nnb
		\end{align}
		Using the fact that  $\alpha\beta\le \frac{\alpha^2}{4} + \beta^2, \forall \alpha,\beta \in \mathbb{R}$ yields
		\begin{align}
		\textstyle \sum_{k=0}^{t} \gamma(k)
		\le 
		\frac{1}{1-\lambda} \sum_{k=0}^{t-1} \big( \frac{\alpha^2(k)}{4} + \beta^2(k) \big). \label{eqBoundSumGamma_t}
		\end{align}
		Moreover, by \eqref{eqBoundBetaSquared}, we have $\beta^2(k) \le n \Phi_k$. Thus 
		\begin{align}
		\textstyle \sum_{0\le k\le t} \beta^2(k)  \le   n \sum_{0\le k\le t} \Phi_k. 
		\label{eqSumBeta2}
		\end{align}
		Using this bound and  \eqref{eqBoundSumGamma_t} in \eqref{eqSumAkXi_Xbar}, we obtain
		\begin{align}
		C_0\!\!\!\!\sum_{0\le k\le t}\!\!\!\!\alpha(k)\| \x_i(k) \!-\! \sB(k)\|
		\!\le\! M_3 \!+\! M_4\!\!\!\!\sum_{0\le k\le t}\!\!\!\!\!\alpha^2(k) \!+\! \frac{D_6}{2}\!\!\!\!\sum_{0\le k\le t}\!\!\!\!\Phi_k \nnb
		\end{align}
		with $C_0 \!=\! \frac{D_6(1-\lambda)}{{2n(Rn+1)C} }$, 
		$M_3 \!=\! \frac{D_6D_{4}\alpha(0)}{2n\lambda C}$, $M_4 \!=\! \frac{D_6}{2n} (  \frac{D_{1}}{\lambda C} +\frac{1}{4})$.
		Combining this inequality with \eqref{eqFbarBound} yields 
		\begin{align}
		&C_0\!\sum_{0\le k\le t}\!\!\alpha(k)\| \x_i(k) \!-\! \sB(k)\| + \!\sum_{0\le k\le t}\!\! \alpha(k)(F(\sB(k)) \!-\! F^*) \nnb\\
		&\!\le (M_1\!+\!M_3) \!+\! (M_2\!+\!M_4)\textstyle \sum_{0\le k\le t}\alpha^2(k), \nnb
		\end{align}
		where $C_1 \!=\! M_1\!+\!M_3, C_2 \!=\! M_2\!+\!M_4$. Dividing both sides by $\sum_{k=0}^t \alpha(k)$ and then using \eqref{eqFEst1} and \eqref{eqFEst2} yields \eqref{eqFConvergence}.
		
		\vspace{3pt}
		{\bf Step (iii): }
		Since  $\tilde{\x}_i(t) \in U$ for $\forall t\ge 0$, $\forall i\in \mathcal{V}$, it then follows from the triangle inequality and Lipschitz continuity of $F$ on $U$ (cf. Assumption \ref{asmp_X}-b) 
		that
		\begin{align} 
		|F(\tilde{\x}_i(t)) - F^*| 
		&\le |F(\tilde{\x}_i(t)) - F(\tilde{\sB}(t))| + F(\tilde{\sB}(t))- F^* \nnb\\
		&\le nL \| \tilde{\x}_i(t) - \tilde{\sB}(t)  \| + F(\tilde{\sB}(t))- F^*.\nnb
		\end{align} 
		Now, by \eqref{eqFConvergence}, both $C_0\| \tilde{\x}_i(t) - \tilde{\sB}(t) \| $ and $F(\tilde{\sB}(t)) - F^* $ are bounded above by $E(t)$. Therefore, \eqref{eqFConvergenceIdenticalConstraints} must hold.%
		\hfill$\blacksquare$
	
	This result demonstrates how the convergence property of the step size sequence implies that of our algorithms. (As a side note, Assumption \ref{asmpStepSize} is not needed here.) 
	The convergence rate analysis now boils down to studying the behavior of $E(t)$; exactly the same task has been carried out thoroughly in the literature for centralized (projected) subgradient methods (see, e.g.,  \cite{Bertsekas99Book,Boyd03Notes,Nesterov04}). 
	Thus, we proceed no further than recalling a few notable results and discussing the constants associated with the convergence rate  in terms of network size and $(1-\lambda)$. 
	\begin{corollary} \label{corComplexity}
		Let the assumptions of Theorem \ref{thmConvergenceRate} hold.
		\begin{itemize}
			\item[(a)] If $\alpha(t)\equiv \alpha$, then $E(t) = C_2\alpha + \frac{C_1}{\alpha t}$. 
			If $\lim_{t\to\infty}\alpha(t) = 0$ and $\sum_{t\ge 0} \alpha(t) = \infty$, then $\lim_{t\to\infty} E(t) = 0$. 
			
			\item[(b)] 
			If $\alpha(t) = O(\frac{1}{\sqrt{t}})$ then $|F(\tilde{\x}_i(t)) \!-\! F^* \!|=\! O(\frac{\ln t}{\sqrt{t}})$.

			\item[(c)]  
			Suppose $\pi^{-1}_{\min} = O(n)$. 
			Then $C_0^{-1} = O(\frac{n^3C}{1-\lambda})$ and 
			$C_2 = O(\frac{n^4C^2}{(1-\lambda)^2} g_{\infty}(D_1, \lambda^2) )$ as $n\to \infty$ and $\lambda \to 1$.  
			Further, if all $X_i, i\in \mathcal{V}$ are  compact, then $C_2 = O(\frac{n^4C^2}{(1-\lambda)^2})$.%
		\end{itemize}
	\end{corollary}
	{\bf Proof.~}  
		%
		%
		We only prove part (c). 
		In view of \eqref{eqEta}, we have $\eta = O(\pi^{-1}_{\min})$. Thus, $\eta = O(n)$. Then, it can be verified that the dominant term is $M_2=  R_2\big(\frac{D_{21}}{1-\lambda}\! + D'_3\big)$ in $C_2$, which is $O\Big(\frac{n^4C^2 R_2}{(1-\lambda)^2}\Big)= O\Big(\frac{n^4C^2}{(1-\lambda)^2}g_{\infty}(D_1, \lambda^2)\Big).$%

		A better estimate can be obtained if we assume further that all $X_i$ are compact. 
		In this case, $\exists D_X>0$ such that $\|\x_{i}(t)-\x^{*} \|^2 \le D_X, \forall i\in \mathcal{V}, \forall t\ge 0$. 
		This enables us to show that \eqref{eqSumUt} in Step (i) of the proof of Theorem~\ref{thmConvergenceRate}  holds with better estimates for $R_1$ and $R_2$. 
		In particular, by using Theorem~\ref{thmMainBound2}, we have for any $t\ge 0$
		\begin{align}
		v_{t+1} &\le v_t  + b_t \textstyle \sum_{i\in \mathcal{V}} \pi_i \|\x_{i}(t) \!-\! \x^{\dagger} \|^2  - u_t + c_t \nnb\\%
		&\le v_t  + D_Xb_t  - u_t + c_t \nnb\\
		&\le v_0  +  \textstyle \sum_{0\le k\le t}  \big(D_Xb_k   - u_k+ c_k\big) \nnb\\
		&\le v_0 + \textstyle \frac{D_X D_1}{1-\lambda^2} +  \sum_{0\le k\le t} (c_k- u_k).\nnb%
		\end{align}
		Here we have used the facts that $\sum_{i\in \mathcal{V}}\pi_i=1$ and $\sum_{k=0}^t b_k = D_1\sum_{k=0}^{t} \lambda^{2k} \le {D_1}/{(1-\lambda^2)}$. 
		As a result, 
		\begin{align}
		\textstyle \sum_{0\le k\le t} u_k \le v_0 +\frac{D_1D_X}{1-\lambda^2} +  \sum_{0\le k\le t} c_k, \quad \forall t\ge 0.\label{eqMainBound_compact1}
		\end{align}
		Thus, we have that \eqref{eqSumUt} still holds but with $R_1 = v_0 +\frac{D_1D_X}{1-\lambda^2}$ and $R_2 = 1$ (compared to $R_2 = D_e/2$ previously). 
		Hence, $C_2$ reduces to $O({n^4C^2}/{(1-\lambda)^2})$.%
		\hfill$\blacksquare$
	
	Note that the rate $O(\frac{\ln t}{\sqrt{t}})$ is also achieved by  distributed subgradient based methods such as Subgradient-Push \cite{Nedic15AC} and Proximal-Gradient \cite[Chap.~3]{Chen12} for unconstrained problems, and Dual Averaging \cite{Duchi12AC} with identical constraints. 
	Here, we have shown that the same rate applies for constrained problems with nonidentical constraint sets, even when $\{ \tilde{\x}_i(t)\}_t$ are infeasible.%
	
	Note from Corollary~\ref{corComplexity}(c) that the spectral gap, defined as $1-|\lambda_2(W)|$, also affects the constant bounds  since $|\lambda_2(W)|\!<\! \lambda \!<\!1$, signifying the importance of strong connectivity. This corollary also suggests the step size be $O(n^{-\sigma})$ (e.g., $\sigma = 2,3$) to lower the order of $E(t)$.

	Finally, we close this section with a remark on the scaling of  $\pi^{-1}_{\min}$. In general, $n \le \pi_{\min}^{-1}$, where equality holds when $W$ is doubly stochastic.  Thus, in Corollary~\ref{corComplexity}(c), we assume  $\pi^{-1}_{\min} = O(n)$ for simplicity of analysis. When $W$ is only row stochastic and $w_{ij} = 1/d_i$ if $(ij)\in \mathcal{E}$ (and $w_{ij} = 0$ otherwise) where $d_i$ is the in-degree of agent $i$, it holds (see \cite{Chung05}) that $\pi_{\min}^{-1} \le n (\max_i d_i)^{\diam(\mathcal{G})}$, where $\diam(\mathcal{G})$ denotes the diameter of $\mathcal{G}$.%
	
	\vspace{-1mm}	
	\section{Numerical example} \label{secNumericalExample}
	\vspace{-1mm}

\begin{figure}[h]
		\centering
		\includegraphics[scale=.4]{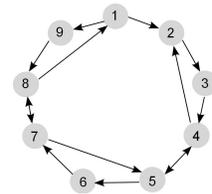} 
		\caption{Communication graph of the network example.} 
		\label{figNetworkExample}
	\end{figure}

	Consider a machine learning problem via the $\ell_1$-norm regularized logistic loss functions 
	\begin{align}
	\textstyle{\min_{\x \in X}}~F(\x) = \textstyle \sum_{i=1}^{r} \ln \big[ 1\!+\! e^{-l_i(\p_i\T \uu + v)} \big] \!+\! \sigma\|\uu\|_1\nnb
	\end{align}
	with variable $\x {=} [\uu\T, v]\T$, $\uu \in \Rm, v \in \mathbb{R}$. Here, $\sigma>0$ is a regularization parameter. 
	The training set consists of $r$ pairs $(\p_i,l_i)$ where $\p_i\in \Rm$ is a feature vector and $l_i \in \{-1,1\}$ is the corresponding label. 
	Suppose that $\x$ satisfies a linear equality constraint: $X = \{ \x\in \mathbb{R}^{m+1}, A_{eq} \x = \bb_{eq} \}$, where $A_{eq} \in \mathbb{R}^{q\times (m+1)}$ and $\bb_{eq} \in \mathbb{R}^q$. 
	In general, when the problem data is distributed or too large to store and/or process on a single machine, employing a network of machines provides a solution. 
	This arises in many applications such as online social network data, wireless sensor networks, and cloud computing. 
	%
	%
	%
	%
	In our example, the problem is to be solved by a network as described in Fig.~\ref{figNetworkExample}.
	
	We assume $r = 500$, $m=50$ and $q=36$, and select $(\p_i,l_i), A_{eq}$ and $\bb_{eq}$ based on normally distributed random numbers. 
	We choose $\sigma = 50$. 
	Suppose the problem data are distributed among the $n$ nodes as follows: each node $i$ stores a partition $\mathcal{P}_i$ of roughly $\frac{r}{n}$ training data and a set of $\frac{q}{n}$ equality constraints, represented by $(A_{eq}^{(i)}, b_{eq}^{(i)})$. 
	Thus, each agent $i\in \mathcal{V}$ is associated with
	\begin{align}
	f_i(\x) &= \textstyle \sum_{j \in \mathcal{P}_i}\! \ln \big[ 1+ e^ {-l_j(\p_j\T \uu + v)}  \big] + \frac{\sigma}{n}\|\uu\|_1\nnb\\
	X_i &= \{ \x\in \mathbb{R}^{m+1}:  A_{eq}^{(i)} \x = \bb_{eq}^{(i)} \}. \nnb
	\end{align}
	Clearly, $X_i$ are unbounded and $X$ has no interior point. We assume that the weight matrix $W=[w_{ij}]$ is such that $w_{ij}= |\mathcal{N}_i|^{-1}$ if $j\in \mathcal{N}_i$ and $w_{ij}=0$ otherwise, where $\mathcal{N}_i$ is the set of node $i$'s in-neighbors (including itself).

	We simulate  Algorithms \ref{algSubgrad1} and \ref{algSubgrad2} using step size $\alpha(t) = {n^{-3}}{(t+1)^{-0.8}}, t\ge 0$ and the usual DPS method (denoted by  DPS-(a)), and its variation DPS-(b) (i.e., the order of the subgradient and consensus steps is reversed) using $\alpha'(t) = {n^{-2}}{(t+1)^{-0.8}}$. 
	Here $\alpha(t)$ and $\alpha'(t)$ differ by a factor of  $n$ for a fair comparison since  subgradients in our algorithms are scaled by $\pi_i^{-1}$ (which equals $n$ if $W$ is doubly stochastic). The initial state $\x_i(0) =\0$ $\forall i\in \mathcal{V}$. Moreover, we use CVX \cite{CVX}  to determine 
	$F^*$ and $\x^*$ by solving the global problem in a centralized fashion. Here, $F^* = 687.67$ and $\|\x^*\| = 1.1341$.

\begin{figure}[t!]
		\hspace{-6mm}\includegraphics[scale = .62]{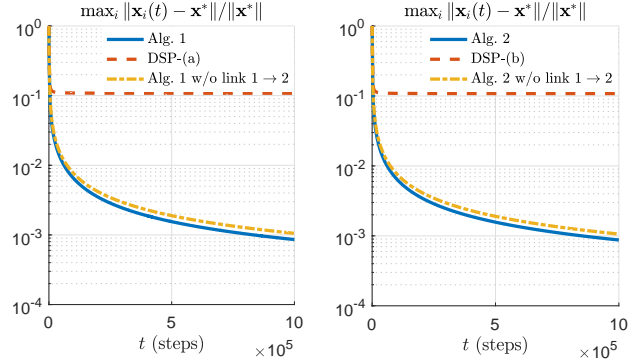}
		\caption{Accuracy of Algorithms \ref{algSubgrad1}, \ref{algSubgrad2}, and DSP methods, where $\|\x^*\| = 1.1341$, DPS-(a): usual  one with consensus step then subgradient, and DPS-(b): subgradient step then consensus.}
		\label{accuracy_compare}
	\end{figure}

	In Fig.~\ref{accuracy_compare}, we show the accuracies of our algorithms and the DPS methods. Clearly, both Algorithms \ref{algSubgrad1} and \ref{algSubgrad2} converge to $\x^*$ and 
	have similar performance, while the DPS methods fail to converge to $\x^*$. We also simulate our algorithms with a new network obtained by deleting link $1\to 2$ in the original network. Here, the algorithms are unchanged except for node $2$ adjusting its incoming link weights. Clearly, convergence is still achieved (since the network is still strongly connected) but slower since the spectral gap is reduced.

	Fig. 3 demonstrates the decay of the objective errors evaluated at the agents' local estimates, namely, $\x_i(t)$ and $\tilde{\x}_i(t)$, when applying the algorithms to both the original and new network. We also show a scaled version of the theoretical  upper bound  in \eqref{eqFConvergenceIdenticalConstraints}, namely $\tilde{E}(t)=(1+nL/C_0)E(t)/F^*/ (3\times 10^{23})$. Although the bound is very loose (due to very rough estimates of $C_i$), the convergence rates of this bound and $|F(\tilde{\x}_i(t)) - F^*|$ are similar, in line with our theoretical findings in Theorem~\ref{thmConvergenceRate}. 
	Again, we  notice that both algorithms have almost the same performance and that the convergence rate is reduced when the network is less connected. 

	\begin{figure}[t!]
		\hspace{-6mm}\includegraphics[scale = .62]{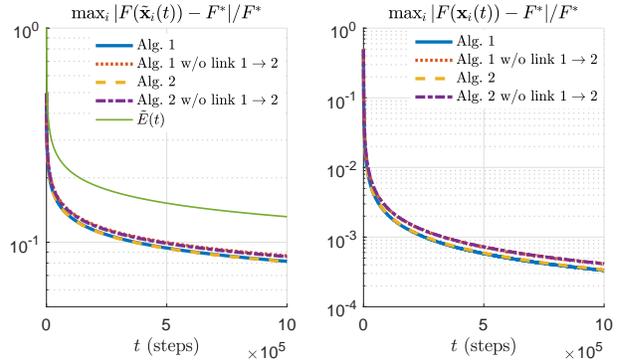}
		\caption{Objective errors evaluated at local estimates $\tilde{\x}_i(t)$ and ${\x}_i(t)$ in Algorithms \ref{algSubgrad1}, \ref{algSubgrad2}, where $F^* = 687.67$ and the green solid line represents $\tilde{E}(t) = (1+nL/C_0)E(t)/F^*/ (3\times 10^{23})$.} 
		\label{objective_error_compare}
	\end{figure}

	\section{Conclusions and discussion} \label{secConclusion}
	In this paper, we have proposed two modified versions of the DPS method that require only a row stochastic weight matrix and studied their convergence and convergence rates. Our analysis also does not invoke a compactness requirement that is usually imposed on the local constraint sets and is able to deal with various scenarios, including  constrained/unconstrained problems, the sets $X_i$ being bounded/unbounded or identical/nonidentical. 
	
	It is important to note the following. 
	First, the idea of using the augmented iteration \eqref{eq_z} to adjust (sub)gradient magnitudes as in \eqref{eq_x} is not only applicable to the distributed projected subgradient methods, but also can be employed to alleviate the condition of the weight matrix being doubly stochastic for some other existing distributed algorithms (using consensus and (sub)gradient steps). For example, we have observed through simulations that the  gradient-based method proposed in \cite{WE10,Shi15} can be modified in the same spirit and still retains fast convergence speed under a suitable constant step size. In \cite{XiMai16}, based on this idea, we proposed a new algorithm that converges linearly under the strong convexity assumption on the cost functions.  
	Second, it is possible to employ other eigenvector estimation schemes in place of \eqref{eq_z} as long as $z_{ii}(t)\to \pi_i$  sufficiently fast (e.g., satisfying \eqref{eqRateBound}). This includes any finite-time computation algorithm, e.g., \cite{Charalambous16}. 
	Third, the convergence analysis developed here can be adapted to either relax the compactness requirement in other projected subgradient based methods (e.g., \cite{Nedic10AC,Lin16distributed}) or impose regular constraints for other subgradient based algorithms (e.g., \cite{Nedic15AC,Olshevsky16}); this holds even when the network is time-varying and possibly with fixed communication delays. 
	Finally, we believe the technique in designing step sizes in \cite{LinRen17arXiv} can enable relaxation of subgradient boundedness, thereby enlarging the addressable class of problems.


	\bibliographystyle{plain}        
	\bibliography{IEEEfull,RefProposal,RefOptim1,RefOptim2}           

\end{document}